\documentclass[10pt]{amsart}

\usepackage{amssymb,latexsym}
\usepackage{amsmath}
\usepackage{graphicx}
\usepackage{amscd}
\usepackage{amsmath}
\usepackage{latexsym}

\usepackage{amsthm}
\usepackage{amsfonts}
\usepackage{latexsym}
\usepackage[all,poly,knot,dvips]{xy}
\usepackage{amsthm,latexsym}
\usepackage{eucal,latexsym}
\usepackage{pstcol} % To use the standard "color" package with PSTricks
\usepackage{pst-node}

\newrgbcolor{mash}{.2 .5 .7}
\newrgbcolor{moosh}{.6 .3 .5}
\definecolor{Pink}{rgb}{1.,0.75,0.8}
\input xy
\newtheorem{theorem}{Theorem}[section]
\newtheorem{proposition}[theorem]{Proposition}
\newtheorem{lemma}[theorem]{Lemma}
\newtheorem{corollary}[theorem]{Corollary}
\newtheorem{definition}[theorem]{Definition}
\newtheorem{example}[theorem]{Example}
\newtheorem{remark}[theorem]{Remark}

\newcommand{\edg}{\text{edg}}
\renewcommand{\vert}{\text{vert}}

\newcommand{\T}{\mathbb{T}}
\renewcommand{\S}{\mathbb{S}}
\newcommand{\R}{\mathbb{R}}
\newcommand{\Z}{\mathbb{Z}}
\renewcommand{\L}{\mathbb{L}}
\renewcommand{\H}{\mathbb{H}}
\newcommand{\D}{\Delta}
\newcommand{\Ker}{\mathop{\rm Ker}}
\numberwithin{equation}{section}
\newcounter{llistadepth}

\newenvironment{manlist}[1]{\addtocounter{llistadepth}{1}
      \edef\llistacontador{llista\romannumeral\the\value{llistadepth}}
      \list{({#1{\llistacontador}})}{\usecounter{\llistacontador}
      \def\makelabel##1{\hss\llap{##1}}
      \itemsep=2pt\parsep=0pt\topsep=3pt plus 1pt minus 1 pt}}{\endlist
      \addtocounter{llistadepth}{-1}}
\newenvironment{alphlist}{\begin{manlist}{\alph}}{\end{manlist}}
\newenvironment{romlist}{\begin{manlist}{\roman}}{\end{manlist}}
\newenvironment{numlist}{\begin{manlist}{\arabic}}{\end{manlist}}

 \psset{unit=0.45}
\title{On String Topology of Three manifolds}
\author{ Hossein Abbaspour}
\address{Centre de Math\' ematique, Ecole Polytechnique, Palaiseau 91128, France.}
 \email{habbaspour@math.polytechnique.fr}

\begin{document}
\maketitle

\begin{abstract}Let $M$ be a closed, oriented and smooth
manifold of dimension $d$. Let $\L M$ be the space of smooth loops
in $M$. In \cite{CS} Chas and Sullivan introduced loop product, a
product of degree $-d$ on the homology of $LM$. In this paper we
show how for three manifolds the ``nontriviality'' of the loop
product relates to the ``hyperbolicity'' of the underlying
manifold. This is an application of the existing powerful tool and
results in three dimensional topology such as the prime
decomposition, torus decomposition, Seifert theorem, torus
theorem.
\end{abstract}
\tableofcontents

\section{Introduction and the Statement of the Main Theorem}
Throughout this manuscript $M$ denotes a connected, oriented
smooth manifold unless otherwise stated. We shall follow the
terminology and notations of \cite{CS}.

 We think of the unit circle, $\S^1$, as the
quotient $\R/\Z$. The \emph{marked point} of the circle is the
image of $0$ in this quotient. A \emph{loop} in $M$ is a
continuous map from $\S^1$ to $M$. The \emph{free loop space} $\L
M$ of $M$ is the space of all loops in $M$. If $f:\S^1 \rightarrow
M$ is a loop in $M$ then the image of the marked point of $\S ^1$
is called the \emph{marked point} of the loop $f$.

We denote the integral homology of $\L M$ by homology  $\H_*(M)$
and is equipped with an \emph{associative product} called
\emph{loop product} (see \cite{CS}) to be denoted by $\bullet$.
The loop product has degree $-d$, where $d$ is the dimension of
$M$ namely,
$$
\bullet : \H _i(M) \otimes \H _j(M)\rightarrow \H_{i+j-d}(M).
$$

\noindent\textbf{Theorem} (Chas-Sullivan \cite{CS}) For a
connected closed oriented smooth manifold $M$, $(\H_{*-d}(M),
\bullet )$ is a graded commutative algebra.

 The inclusion of constant loops induces a map
$$
i:H_*(M)\rightarrow \H _*(M),
$$
and $p:\H _*(M)\rightarrow H_*(M)$ is the map induced by
projection on the marked points. Observe that
$$
p\circ i=id_{H_*(M)},
$$
 hence there is a canonical decomposition
\begin{center}
$\H _*(M)=i(H_*(M)) \oplus A_M, $\\
$A_M=\Ker p$ and $ p_{A_M}=id_{\H _* (M)}- p.$
\end{center}
\begin{proposition}(See \cite{CS})
\label{proposition1}
 $i:H_*(M)\rightarrow \H _*(M)$ is a map of graded algebras where  the multiplication on $H_*(M)$ is
the intersection product denoted by $\wedge$.
 Therefore, one can regard $(\H_* (M),\bullet)$ as an extension of $(H_* (M), \wedge)$.
\end{proposition}
\begin{definition}
\emph{We shall say that an oriented closed manifold $M$ has
\emph{nontrivial extended loop products} if the restriction of the
loop product to $A_M$ is nontrivial, namely $$
\bullet|_{A_M\otimes A_M}\neq0.$$ Otherwise we shall say that $M$
has \emph{trivial extended loop products}.}
\end{definition}

\begin{definition}\emph{
Let $M$ be a closed $3$-manifold. We say that $M$ is
\emph{algebraically hyperbolic} if it is a $K(\pi,1)$ (aspherical)
and its fundamental group has no rank $2$ abelian subgroup.}
\end{definition}
Note that any finite cover of a algebraically hyperbolic manifold
is also algebraically hyperbolic.\\ The following is the main
theorem of this paper:
\begin{theorem}
\label{theorem1} Let $M$ be an oriented closed $3$-manifold.
If $M$ is not algebraically hyperbolic then $M$ or a double cover of $M$ has nontrivial extended loop products.\\
 If $M$ is algebraically hyperbolic then $M$ and all its
 finite covers have trivial extended loop products.
\end{theorem}
  From now on we identify $H_*(M)$
with its image under $i$. The following lemma will be useful to
us.

\begin{lemma}:\label{lem-compoenents}
\begin{romlist}
\item  There is a one-to-one correspondence between the connected
components of $\L M$ and the conjugacy classes of $\pi_1(M)$.
\item Let $\alpha \in \pi_1(M)$ and let $f \in \L M$ be a loop
representing $\alpha$. Then there exists a short exact sequence
\begin{equation}
\begin{CD}
\pi_2(M)\cong \pi_1((\Omega M)_{\alpha},f) @>S>>
\pi_1((\mathbb{L}M)_{[\alpha]}, f) @>R>> \pi_{\alpha} @>>> 0
\end{CD}
\end{equation}
where $\pi_{\alpha}$ is the centralizer of $\alpha$ in $\pi_1(M)$.
\end{romlist}
\end{lemma}

 If
$[\alpha]\in \hat{\pi}_1(M)$ is a conjugacy class then $(\L
M)_{[\alpha]}$ denotes the corresponding connected component of
$\L M$ and $\H _*(M)_{\alpha}$ its integral homology.

If $\alpha$ and $\beta$ are two loops with the same marked points,
their \emph{composition}, $\alpha \beta$, is also a loop and
 is defined
 as follows:
\begin{equation*}
\alpha\beta=
\begin{cases}
 \alpha(2t) &  0\leq t\leq 1/2\\
 \beta(2t-1) &  1/2\leq t \leq1\\
\end{cases}
\end{equation*}

We give an informal description of loop product: given two
homology classes in $\H _*(M)$, choose a chain representative for
each one of them and consider the corresponding set of marked
points. They form two chains in $M$. If these two chains are not
in transversal position, replace the representatives with
appropriate ones so that their chain of marked points intersect
transversally. At each intersection point, compose the
corresponding loops of each chain of loops. This process produces
a new chain of loops. In \cite{CS} it is proved that this
operation passes to homology.

\begin{remark}\emph{If $M$ is a closed manifold then the algebra $(\H
_*(M),\bullet)$ has a \emph{unit}, $\mu_M$, described as follows.
Since each point in $M$ can be regarded as a constant loop, a
chain in $M$ gives rise to a chain in $\L M$. This map passes to
homology. The image of the fundamental class of $M$ under this map
is the unit of $(\H _*(M),\bullet)$.}
\end{remark}

 $\H_* (M)$ is also equipped
with a self map of degree $1$, denoted $\D$,  which can be
described as follows. Consider the natural action of $\S^1$ on $\L
M$: $f \in \L M$ and $r\in \S^1=\R/\Z$,
$$
r.f: \S^1 \rightarrow M
$$
$$
 (rf)(t)=f(r+t)
$$
This action induces a map
$$
\D : \H _*(M) \rightarrow \H _{*+1}(M)
$$
$\D$ has the property that $\D ^2=0$ (see \cite{CS}).

\begin{example}
\label{example1}(see \cite{CS}) \emph{Let $M$ be a closed
hyperbolic oriented
 $3$-manifold. As $M$ is aspherical then by
Lemma \ref{lem-compoenents}, each component of $(\L M)_{[\alpha]}$
is a $K(\pi_{\alpha},1)$ where $\pi_{\alpha}$ is the centralizer
of $\alpha$ in $\pi_1(M)$. If $\alpha\neq 1 \in \pi_1(M)$ then it
can be regarded as a hyperbolic element of the discrete subgroup $
\pi_1(M)$ of $PSL(2,\mathbb{C})$. So its centralizer is the cyclic
group generated by $\alpha$. Hence $(\L M)_{[\alpha]}$ is a
$K(\Z,1)$ and $$A_M \simeq H_*(\underset{[\alpha]\neq [1] \in
\hat{\pi}_1(M)}{\coprod} K(\Z,1)).$$ Therefore $A_M$ has
homological dimension $1$ which implies that the loop product is
identically zero on $A_M\otimes A_M$ as a consequence of the
absence of homology class of sufficiently large dimension. So
closed oriented hyperbolic manifolds have trivial extended loop
product. }
\end{example}
\section{3-Manifolds with Finite Fundamental Group}
 In this section we show that $3$-manifolds with
finite fundamental group have nontrivial extended loop products.

The following construction is valid for any oriented closed
manifold $M$. Let $p$ be the based point for $M$. Consider the
$0$-chain in $\L M$ that consists of the constant loop at $p$. It
represents a homology class $\epsilon \in \H _0(M)$. Given a chain
$c\in \L M$ whose set of marked points is transversal to $p$, then
the loops in $c$ whose marked points are $p$ form a chain in the
based loop space of $\Omega M$. This induces a map of degree $-d$
$$
\cap:\H _*(M)\rightarrow H _{*-d}(\Omega M).
$$
On the other hand $H_*(\Omega M)$ is equipped with Pontryagin
product.

\begin{lemma}(see \cite{CS})
\label{lemma6}
 $\cap:(\mathbb{H}(M),\bullet )\rightarrow
(H_{*-d}(\Omega M),\times) $ is an algebra map where $\times$ is
the Pontryagin multiplication. In particular if $M$ is a Lie group
then $\cap$ is surjective.
\end{lemma}
\begin{proposition}
\label{S3}
 $S^3$ has nontrivial extended loop products.
\end{proposition}

\begin{proof}
Since  $S^3$ is a Lie group, by Lemma \ref{lemma6} $\cap$ is
surjective. $(H_*(\Omega S^3), \times) $ is a polynomial algebra
with one generator of degree $2$.  Let  $x_1\in H_p(\Omega S^3)$
and $x_2 \in H_q(\Omega S^3)$ be non-zero elements such that $p,q
>3$. Let $a_i\in \cap^{-1}(x_i)$ for $i=1,2$. Since $H_i(S^3)=0$ for $i>3$ we have $a_i \in
A_{S^3}$ . By Lemma \ref{lemma6} $\cap(a_1\bullet
a_2)=\cap(a_1)\times \cap_p(a_2)=x_1 \times x_2\neq 0$, therefore
$a_1 \bullet a_2\neq 0$
\end{proof}

\begin{lemma}
\label{lemma7} Let $M$ and $N$ be two closed oriented
$n$-manifolds and $p\in N$.
 For $f:M\rightarrow N$ which is transversal to $p$ and
$f^{-1}(p)$ consists of only one point, the following diagram
commutes :
\begin{equation}
\label{commu-homot-S3}
\begin{CD}
\mathbb{H}_*(M) @>{\cap}>> H_*(\Omega M) \\
@V{f_{\mathbb{L}}}VV  @V{f_{\Omega}}VV\\
 \mathbb{H}_*(N) @>{\cap}>> H_*(\Omega N)
\end{CD}
\end{equation}
$f_{\Omega}$ and $f_{\mathbb{L}}$ are the induced maps by $f$ on
the homologies of the based loop space and the free loop space.
\end{lemma}
\begin{proof}
We choose $p$ and $f^{-1}(p)$ respectively as the base point of
$N$ and $M$. We show that the above diagram is commutative even at
the chain level. Let $k:\Delta ^m \rightarrow \L M$ be a
\emph{m}-simplex in $\L M$. So for each $x\in\Delta ^m$
$$
k(x):\S ^1 \rightarrow M.
$$

\noindent If $k_0: \Delta ^m \rightarrow M$,
$$
k_0(x)=k(x)(0),
$$
is transversal to $f^{-1}(p)$ then consider the set $\Lambda
\subset \Delta^m$,

$$
\Lambda=\{x\in \Delta ^m|k(x)(0)=f^{-1}(p)\}.
$$

 \noindent $k|_{\Lambda}$ is a chain in $\Omega_{f^{-1}(p)} M$
and $\cap(k)=k|_{\Lambda}$. After composing with $f$, $f\circ
k|_{\Lambda}$ is a chain in $\Omega_p N$ and
$$
(f_{\Omega}\circ\cap)(k)=f\circ k|_{\Lambda}.
$$

On the other hand $f_{\L }(k)=f \circ k$ is a chain in $\L N$.
Since $k_0$ is transversal to $f^{-1}(p)$ and $f$ is transversal
to $p$ thus $f\circ k_0$ is transversal to $p$. Since $f^{-1}$
consists of only one point in $M$ then
$$
\{x\in\Delta ^m| ((f\circ k)(x))(0)=p\}=\{x\in \Delta
^m|k(x)(0)=f^{-1}(p)\}=\Lambda,
$$
therefore,
$$
(\cap \circ f_{\L})(k)= f\circ k|_{\Lambda}=(f_{\Omega}\circ \cap)
(k).
$$

\end{proof}

\begin{lemma}
\label{lemma-homotopy-equi} Let $M$ be a 3-manifold which is a
homotopy $3$-sphere and $p_1,p_2,...,p_r$ are distinct points in
$M$. Then there is a homotopy equivalence $f:S^3 \rightarrow  M$
which is transversal to all $p_i$'s and $f^{-1}(p_i)$ consists of
only one point for each $i$.
\end{lemma}
\begin{proof}
 Consider a closed ball $D\subset S^3$. Let $f$ be a diffeomorphism from $D$ to a closed ball in
$D'\subset M$ where $p_i\in D'$ for all $i$.

$f|_{\partial D}:\partial D\rightarrow\partial D'$ can be extended
to  $f:S^3\setminus Int D' \rightarrow M \setminus Int D$ as all
the obstructions which are cohomology classes in
$$
H^{q+1}(S^3\setminus Int D',\partial D',\pi_q(M \setminus Int D))
$$
 vanish since $M \setminus Int D'$ is contractible\footnote{This can be proved using Mayer-Vietrois exact sequence}.

Therefore we have a map that $f:S^3\rightarrow M$ which sends $D'$
to $D$ diffeomorphically and sends the complement of $D'$ to the
complement of $D$ and in particular it is transversal to $p_i$ and
$f^{-1}(p_i)$ consists of only one point in $S^3$.
 To see that $f$ is a homotopy equivalent it is
enough to observe that it has degree one.
\end{proof}
\begin{proposition}
\label{proposition3}
 A closed simply connected  $3$-manifold $M$, has nontrivial extended loop products.
\end{proposition}
\begin{proof}
Let $p\in M$ and $f:S^3\rightarrow M$ be the homotopy equivalence
provided by Lemma \ref{lemma-homotopy-equi} for $r=1$.

 $f$ induces a homotopy equivalences $f_{\Omega}:\Omega S^3 \rightarrow \Omega M$ and
$f_{\mathbb{L}}:\mathbb{L}S^3  \rightarrow \mathbb{L}M$. We have
the following commutative diagram by Lemma \ref{lemma7}:
\begin{equation}
\label{commu-general}
\begin{CD}
\mathbb{H}_*(S^3) @>{\cap_{S^3}}>> H_*(\Omega_{f^{-1}(p)} S^3)\\
@V{f_{\mathbb{L}}}VV  @V{f_{\Omega}}VV \\
\mathbb{H}_*(M) @>{\cap_{M}}>> H_*(\Omega _pM)\\
\end{CD}
\end{equation}
Since $f_{\Omega}$ and $f_{\mathbb{L}}$ are isomorphisms and
$\cap_{S^3}$ is surjective, therefore $\cap_{M}$ is also
surjective.
 Let $x_1\neq 0 \in H_m(\Omega M)$ \footnote{ $H_*(\Omega
M)\cong H_*(\Omega S^3)\cong\Z[x]$
 and degree $x$=2}and $x_2\neq 0  \in
H_n(\Omega M)$
 where $n,m>3$ and $a_i\in \cap^{-1}(x_i)$, $i=1,2$.

Then, $a_i\in A_M$ since $H_k(M)=0$ for all $k>3$. By
commutativity of diagram (\ref{commu-general}),
$$
f(a_1\bullet a_2)=f(a_1)\times f(a_2)=x_1\times x_2 \neq 0 .
 $$
Hence $a_1\bullet a_2\neq 0.$

\end{proof}
\begin{proposition}
\label{finitegroup} If $M$ is a closed oriented $3$-manifold with
finite fundamental group then $M$ has nontrivial extended loop
products.
\end{proposition}

\begin{proof}
Let $\tilde{M}$ be the universal cover of $M$ and $q$ be the
covering map and  $r=degree(q)=|\pi_1(M)|< \infty$. Choose $p\in
M$ and let $\{p_1,p_2,...,p_r\}\in q^{-1}(p)$ where $p_i$'s are
distinct.

Since $\tilde{M}$ is homotopy equivalent to $S^3$, by Lemma
\ref{lemma-homotopy-equi} there is a homotopy equivalence $f:S^3
\rightarrow \tilde{M}$ which is transversal at each $p_i$ and
$f^{-1}(p_i)$ consists of only one point in $S^3$ in $M$. Let
$m_i=f^{-1}(p_i)$, $i=1,2,...,r$. Hence $(q\circ
f)^{-1}(p)=\{m_1,m_2,...,m_r\}$ and $q\circ f$ is transversal to
$p$.

\noindent$f$ composed by the covering map $q:\tilde{M}\rightarrow
M$ induces the maps
$$(q\circ
f)_{\mathbb{L}}:\mathbb{H}_*(S^3)\rightarrow \mathbb{H}_*(M)$$
\begin{center} and \end{center} $$(q\circ
f)_{\Omega}:H_*(\Omega_{m_1} S^3)\rightarrow H_*(\Omega_p M)$$

\noindent where the latter is an isomorphism.

As we know, generally $\Omega M$ and $\mathbb{L}M$ are not
connected. So we should emphasize that the images of $(q\circ
f)_{\mathbb{L}}$ and $(q \circ f)_{\Omega}$ are homologies of the
connected components corresponding to the trivial loop in $M$ and
in the diagram below, the sub-index $0$ indicates that.

 We claim that
\begin{equation}
\label{comm-finite}
\begin{CD}
\mathbb{H}_*(S^3) @>{r.\cap}_{S^3}>> H_*(\Omega_{m_1} S^3) \\
@V{(f\circ q)_{\mathbb{L}}}VV  @V{(f\circ q)_{\Omega}}VV\\
(\mathbb{H}_*(M))_0 @>{\cap}_M>> H_*((\Omega_p M)_0)
\end{CD}
\end{equation}
is commutative. To prove this note that $S^3$ has a group
structure hence
$$\mathbb{L}S^3\overset{homeo.}{\simeq}S^3 \times \Omega S^3,$$
and
\begin{equation}
\label{eq-S3}
 \mathbb{H}_*(S^3)\simeq H_*(S^3)\otimes
H_*(\Omega S^3).
\end{equation}
There are two types of homology classes in $\mathbb{H}_*(S^3)$.
\begin{numlist}
\item The classes that under the isomorphism (\ref{eq-S3})
correspond to homology classes in $H_0(S^3)\otimes H_*(\Omega
S^3)$.  These classes can be represented by cycles in
$\mathbb{L}S^3$ whose sets of marked points consist of a single
point in $S^3$.

 \noindent Type (1) homology classes are mapped to
$0$ by $\cap_{S^3}$, since we can choose a chain representative
whose set of marked point, which consists of one points, is
different from all $m_i$'s. Therefore they are mapped to $0$ by
$(f\circ q)_{\Omega}\circ({r.\cap}_{S^3})$

 \noindent Indeed their image under $(f\circ q)_{\mathbb{L}}$
have sets of marked point consist of only one point, $p$. Thus,
they are mapped to zero under ${\cap}_M \circ(f\circ
q)_{\mathbb{L}}$. Hence the diagram commutes in this case

\item The classes that under isomorphism (\ref{eq-S3}) correspond
to homology classes in $H_3(S^3) \otimes H_*(\Omega S^3)$. They
can be represented by a cycle in $\mathbb{L}S^3$ whose set of
marked points is all of $S^3$.

\noindent Let $\theta=a\otimes b \in \H_3(S^3)\simeq  H_3(S^3)
\otimes H_*(\Omega S^3)$. Then
$$
\cap_{S^3}(\theta)= b
$$
so
$$
(f\circ q)_{\Omega}\circ(r\cap_{S^3}\theta)=r(f \circ q)(b)
$$
and this is exactly $\cap_M\circ(f\circ q)_{\Omega}(b)$. Note that
under $f\circ q$, $r$ copies of $b$ which are based at each one of
$m_i$'s, come together.
\end{numlist}

Now consider $a_1,a_2\in \H _*(S^3)$ constructed in the proof of
Proposition \ref{S3}. From commutativity of diagram
(\ref{comm-finite}) and the fact that $(q\circ f)_{\Omega }$ is an
isomorphism it follows that
 $$(q\circ f)_{\L}(a_1) \bullet
(q\circ f)_{\L}(a_2)\neq 0$$
\end{proof}
\section{3-Manifolds with Non Separating 2-Sphere or 2-Torus}
\label{nonseparating-section}
\begin{proposition}
\label{S1S2}
 $S^1\times S^2$ has nontrivial extended loop products.
\end{proposition}
\begin{proof}
Consider $a \in S^1$, $b,p \in S^2$ and $p\neq b$ and let $(a,p)$
be the base point of $S^1\times S^2$. The maps:
$$
\alpha : S^1 \rightarrow S^1\times S^2 ,\qquad \alpha(x)=(x,b),
$$
gives rise to a loop whose marked point is $(a,b)$.

$\alpha$ represents a 0-homology class $\mathbb{H}_0(M)$. The loop
$\gamma \alpha \gamma^{-1}$ represents an element $g_1\in
\pi_1(M)$, where $\gamma$ is a path
 connecting $(a,b)$ and $(a,p)$. The free homotopy type of $\alpha$ is $[g_1]\in \hat{\pi}_1(M)$.

Consider the map

$$
\beta: S^2 \rightarrow S^1\times S^2 ,\qquad \beta(y)=(a,y).
$$

 Also, $\beta$ can be regarded as an element of $
\pi_1(\Omega (S^1\times S^2))\simeq \pi_2(S^1\times S^2)$ or in
other words we think of $\beta$ as a $1$-dimensional family of
loops. The free homotopy type of the loops of $\beta$ is the one
of the trivial loop. We compose each loop of $\beta$ with a fixed
loop $\gamma'$ whose marked point is $(a,p)$ and it represents a
nontrivial element $g_2 \in \pi_1(S^2\times S^1)$, $g_2\neq g_1$.
We continue to denote the result by $\beta$.

\begin{center}
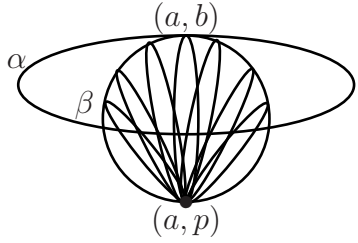
\begin{figure}[h]
\begin{pspicture}(0,-1)(15,5)
 \pscircle[linewidth=1pt](7,2.5){2.5}
\psellipse[linewidth=1.pt](7,2.5)(0.4,2.5)
%\psellipse[linewidth=0.5pt](2,2)(0.3,1)
%\rput{30}(0,0){\psellipse[linewidth=0.5pt](2,2)(0.3,1)}
\rput[b]{26.51}(0.8,-3.4){\psellipse[linewidth=1.pt](7,2.5)(0.3,2.23)}
\rput[b]{-26.51}(0.6,2.9){\psellipse[linewidth=1.pt](7,2.5)(0.3,2.23)}
\rput[b]{11.76}(0.06,-1.49){\psellipse[linewidth=1.pt](7,2.5)(0.35,2.45)}
\rput[b]{-11.76}(0.06,1.44){\psellipse[linewidth=1.pt](7,2.5)(0.35,2.45)}
\rput[b]{38.659}(1.9,-4.8){\psellipse[linewidth=1.pt](7,2.5)(0.3,1.920)}
\rput[b]{-38.659}(1.2,3.9){\psellipse[linewidth=1.pt](7,2.5)(0.3,1.920)}
\psellipse[linewidth=1.pt](7,3.5)(5,1.5)
\rput[b]{0}(2,4){\Large$\alpha$}
\rput[b]{0}(7,-0.2){\Large$\bullet$}
\rput[b]{0}(4,2.5){\Large$\beta$} \rput[b]{0}(7,5){\Large$(a,b)$}
\rput[b]{0}(7,-1){\Large $(a,p)$}
\end{pspicture}
\caption{A non separating $2$-sphere}\label{figS2S1}
\end{figure}
\end{center}
The image of $\beta$ under the Hurewicz map
$$
\pi_2(\Omega (S^2\times S^1)) \rightarrow H_2(\Omega( S^2\times
S^1))
$$
and then followed by the induced map by inclusion
$$
H_2(\Omega( S^2\times S^1))\rightarrow \H_2(S^2\times S^1)
$$
 \noindent gives rise to an element in
$\mathbb{H}_2(S^1\times S^2)$. Again we continue to denote this
element by $\beta$.

\noindent We claim  that $$p_{A_{S^1 \times S^2}} (\Delta
\alpha)\bullet p_{A_{S^1\times S^2}}( \Delta \beta ) \neq 0 \in
\mathbb{H}_0(S^1\times S^2) $$
 and this proves that $S^1 \times S^2$ has nontrivial extended loop products.

\noindent Note that
\[ p_{A_{S^1 \times S^2}}( \Delta
\alpha)\bullet p_{A_{S^1\times S^2}}( \Delta \beta )=(\Delta
\alpha - p_{S^1\times S^2}(\Delta \alpha))\bullet (\Delta \beta -
p_{S^1\times S^2}(\Delta \beta)).
\]

In expanding this expression, we get four terms, each one is the
loop product of two homology classes, a $1$-homology class and a
$2$-homology class. The set of marked points of the chain
representatives of the homology classes involved in each one of
these four terms intersect exactly at one point, namely $(a,b)$.
Hence $ p_{A_{S^1 \times S^2}}( \Delta \alpha)\bullet
p_{A_{S^1\times S^2}}( \Delta \beta ) $ is equal to the sum of $4$
homology classes in $\mathbb{H}_0(S^1\times S^2)$, each
represented by a loop and a $\pm$ sign.

 The free homotopy types of these loops are
$[g_1g_2],[g_1],[g_2], 1$ which are different. Therefore the loop
product is nontrivial.
\end{proof}
\begin{corollary}\label{nonseparating2}
 An oriented closed $3$-manifold with a non separating $2$-sphere
has nontrivial extended loop products.
\end{corollary}
\begin{proof}
The proof is similar to the one of Proposition \ref{S1S2}. The
only fact we used in the proof of Proposition \ref{S1S2} was that
there was a non separating $2$-sphere (a $2$-sphere that meets
some closed simple curve only once). Every oriented closed
$3$-manifold with non separating $2$-sphere is indeed a connected
sum of $S^2\times S^1$ with another $3$-manifold (see \cite{He}).
So $\pi_1(M)$
 is the free product of $\pi_1(S^2\times S^1)$ and the fundamental group of another $3$-manifold. Therefore
the choices $g_1,g_2\neq 1$, $g_2\neq g_1$ of elements of
$\pi_1(M)$, as in the proof of Proposition \ref{S1S2}, are also
possible in this case.
\end{proof}
\begin{proposition}
\label{selfglu-torus}
 If a closed oriented $3$-manifold $M$ has a non separating two sided\footnote{The normal bundle is trivial.}
 incompressible\footnote{$\pi_1$-injective.} torus,
 then $M$ has nontrivial extended loop products.
\end{proposition}
\begin{proof}
Let $\varphi:S^1\times S^1 \rightarrow M$ be a non separating
torus. $H$ gives rise to a homology class $H$ in
$\mathbb{H}_1(M)$, where the marked points are $\varphi(t,1),\quad
t\in S^1 $ and the loop passing through $\varphi(t,1)$ is
 \begin{equation*}
 \varphi_t:\S^1 \rightarrow M
  \end{equation*}
   \begin{equation*}
 \varphi_t(s)=\varphi(t,s).
 \end{equation*}

 The loop $\varphi_1$ represents an element $h\in\pi_1(M)$ which is
 nontrivial as the torus is incompressible.
\bigskip
\begin{center}
\begin{figure}[h]
\begin{pspicture}(-1,-3)(15,6)
%\psgrid(0,-6)(15,6)
\psellipse[linewidth=1.pt](6.9,2.5)(1.8,3.3)
\psellipse[linewidth=1.pt](6.9,2.5)(1.2,3.0)
\psellipse[linewidth=1.pt](6.9,2.5)(1.4,3.2)
\psellipse[linewidth=1.pt](6.9,2.5)(1.6,3.3)
\pscurve[linewidth=1.pt](7,4.2)(6.8,3.8)(6.4,2.5)(6.7,1.2)(7,0.8)
\pscurve[linewidth=1.pt](6.9,4)(7.2,2.5)(6.8,1.1)
\pscurve[linewidth=1.pt](1,5.7)(2,5.5)(7,5.8)(9,5.3)(11,5.5)(12,5,3)
\rput[b]{0}(0,-6.6){\pscurve[linewidth=1.pt](1,5.7)(2,5.5)(7,5.8)(9,5.3)(11,5.5)(12,5,3)}
\psccurve[linewidth=1.pt](8,3)(4,3)(0.4,-0.4)(7,-3)(13.4,-0.4)(8,3)

\rput[b]{0}(5.8,3.2){\small$\bullet$} \rput[b]{0}(9,1){\Large $M$}
\rput[b]{0}(11,2.3){\Large $l$}
\psellipse[linewidth=1.pt](5.75,2.6)(0.65,0.2)
\rput[b]{0}(2,3.5){\small  intersection point}
\rput[b]{0}(3,0.6){\small marked points }\rput[b]{0}(3,0.1){of
$H$} \psline{->}(2.7,3.5)(5.7,3.3) \psline{->}(3,1.4)(5.5,2.4)
\end{pspicture}
\caption{A non separating torus}\label{fig-nonsep-torus}
\end{figure}
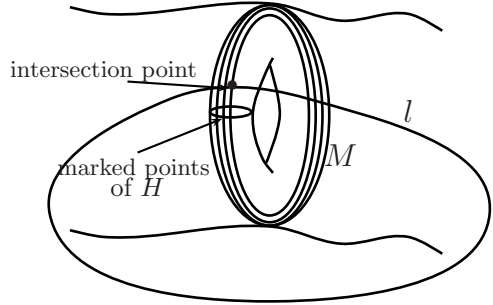
\end{center}
 Let $\lambda:S^1 \rightarrow M $ be a loop in $M$ which intersects the
non separating torus transversally at exactly one point $m$ (see
Figure \ref{fig-nonsep-torus}). We can assume that the
intersection point $m$ is the marked point of $\lambda$ and we
consider this point as the base point of $M$. Without loss of
generality we can assume $m=\varphi(1,1)$.

$\lambda$ represents a nontrivial element $l \in\pi_1(M)$ since
the homology class it represents is nontrivial.

 Also, $\lambda$ gives rise to a homology class $L\in \H_0(M)$. The free homotopy type of the
loops of $L$ is $[l]\in \hat{\pi}_1(M)$.

 We claim that
\begin{equation*}p_{A_M}(\Delta L)\bullet p_{A_M}(\Delta H) \neq 0 \end{equation*} which
proves the assertion.

 In the expansion of $p_{A_M}(\Delta L)\bullet p_{A_M}(\Delta H)
\in \mathbb{H}_0(M)$  we get four terms, each represented by a
single loop and $+$ or $-$ sign. The loops have free homotopy type
$[l.h]$, $[l]$, $[h]$ and $[1]$. If we show that three out of four
terms are different then there cannot be a complete cancellation,
thus the product is not trivial.

Assume $[l]$ and $[h]$ are the same. Then the $1$-homologies
represented by $\lambda$ and $\varphi_1$ are the same. Note that
the $1$-homology  class represented by $\lambda$, intersects the
$2$-homology represented by the torus $\varphi(.,.)$ nontrivially,
while the $1$-homology class represented by $\varphi_1$ intersects
the same 2-homology class trivially \footnote{Since the torus is
two sided (it normal bundle is trivial) and we can push
$\varphi_1$ off the torus.}. This is a contradiction, hence $[l]$
and $[h]$ are different. Similar argument shows that $[l]\neq
[1]$. To see that $[h]\neq [1]$ note that the non separating torus
does not bound in $M$.
\end{proof}
\section{Closed Seifert Manifolds}
\label{section-closed-seifert}
 In this section we prove that closed Seifert manifolds
have nontrivial extended loop products. Before we state and prove
the result we recall the definition of Seifert manifolds.

A $3$-manifold is  called \emph{Seifert manifold} if it admits a
foliation by circles $S^1$ in such a way that each leaf has a
foliated tubular neighborhood $D^2\times S^1$ which is $D^2\times
S^1 $ either with its trivial foliation or its quotient by a
standard action of a finite cyclic group $\mathbb{Z}_p$ (acting on
$D^2$ factor). The fibers (leaves) which have a trivially foliated
$ D^2\times S^1$ neighborhood are called \emph{normal fibers} and
the others are called \emph{singular fibers} and $p$ in above is
its \emph{multiplicity}. The space of the fibers forms a surface
and it is called the \emph{base surface}.

\textbf{Notation:} Let $M$ be an oriented Seifert manifold with
$p$ exceptional fibers and $k$ boundary components. If $S$, the
base surface of $M$, is orientable, then it has genus $g\in \Z^+$
otherwise $S$ is non orientable and $g$ is the number of cross
caps in $S$.
$\pi_1(M)$ has the presentation (see \cite{J}):\\
(1) If $S$ is orientable and of genus $g$,
\begin{equation}
\label{seifert-oriented}
\begin{split}
\pi_1(M)=\langle a_1,b_1,...a_g,b_g,c_1,...c_p,d_1,d_2...d_k,h &|a_i h a_i^{-1}= h,\quad  b_ihb_i^{-1}=h,\\
    & \quad c_ihc_i^{-1}=h, \quad d_ihd_i^{-1}=h, \quad\\ & \quad c_i^{\alpha_i}h^{\beta_i}=1\\
 &\quad 1=\prod^g_{i=1} [a_i,b_i] \prod^p_{i=1} c_i  \prod^k_{i=1} d_i\rangle,
\end{split}
\end{equation}
(2) If $S$ is nonorientable and has $g$ cross caps.

\begin{equation}
\label{seifert-nonoriented}
\begin{split}
\pi_1(M)=\langle a_1,...a_g,c_1,...c_p,d_1,d_2...d_k,h &|a_i h a_i^{-1}= h^{-1},\\
    & \quad c_ihc_i^{-1}=h^{\delta_i}, \quad d_ihd_i^{-1}=h, \quad \\
    & \quad c_i^{\alpha_i}h^{\beta_i}=1\\
 &1=\prod^g_{i=1} a_i^2 \prod^p_{i=1} c_i \prod_{i=1}^kd_i\rangle
\end{split}
\end{equation}
where  $0 < \beta_i < \alpha_i$ are integers and $\delta_i=\pm1$.
$h$ corresponds to the normal fiber. Each $c_i$ can be represented
by a loop on the base surface going around the singular fiber
once. The $d_i$'s correspond to the boundary components.
$\alpha_i$ is the multiplicity of the singular fiber corresponding
to $c_i$ and we say that it has type $(\alpha_i,\beta_i)$. Note
that $\langle h \rangle$ is a normal subgroup of $\pi_1(M)$ and it
is central if the $S$ is orientable.

\begin{remark}
\label{rem-oriented-seifert}\emph{An orientable Seifert manifold
$M$ may not be oriented as a fibration but it always has a double
cover which is Seifert and is orientable as a fibration or in
other words its base surface is orientable. This double cover
$\tilde{M}$ can be described  as following:
$$
\tilde{M}=\{(m,o)| m\in M \quad \& \quad o \text{ an orientation
for the fiber passing through } m \}
$$
The covering map $\tilde{M}\overset{q}{\rightarrow} M$ is
$$
q(m,o)=m.$$ The Seifert fibration of $\tilde{M}$ is obtained by
pulling back the one of $M$.
 If $M$ has boundary the then each boundary component, a torus, gives
rise to $2$ boundary components each of which is a torus. If $M$
has $g$ cross caps then $\tilde{M}$ has genus $g-1$.}
\end{remark}

\begin{proposition}
\label{seifert-closed}
 If $M$ is a closed oriented Seifert $3$-manifold then $M$ or a double cover of $M$
 has nontrivial extended loop products.
\end{proposition}

\begin{proof}
By Remark \ref{rem-oriented-seifert} we may assume that the base
surface of $M$ is orientable.
 Also we assume $\pi_1(M)\ncong \Z_2$, since Corollary \ref{finitegroup} deals with this case. Let $h$ be the generator
of $\pi_1(M)$ corresponding to the normal fiber according the
presentation (\ref{seifert-oriented}).

 There is a natural $3$-homology class in $\H_3(M)$ which has a representative
 whose set of marked points is exactly $M$ and the homotopy type of the loops is $[h]$. The loop
associated with a point in $M$ that is on a normal fiber, is the
fiber passing through the point. The loop passing through a point
on a singular fiber as a map is a multiple of the singular fiber.
We denote this homology class by $\Theta_M$. Loops of $\Theta_M$
have the free homotopy type $[h]$.
 We have $p_{A_M}(\Theta_M)= \Theta_M-p_M(\Theta_M)=\Theta_M-\mu_M$
where $\mu_M$ is the unit of $\H_*(M)$.

Consider $l$ a loop representing an element $\alpha \in \pi_1(M)$,
$\alpha\neq h,1$. This is possible since we are assuming
$\pi_1(M)\neq \Z_2$. Since $h$ is central, the free homotopy class
$[\alpha]$ is different from $[h]$ and $[1]$.

$\alpha$ represents a $0$-homology class $\hat{\alpha}$ in
$\mathbb{H}_0(M)$.
  We have
\begin{equation*}
p_{A_M}(\hat{\alpha})= \hat{\alpha}-
p(\hat{\alpha})\end{equation*}
 where $p(\hat{\alpha})$ is just the
constant loop at the marked point of $l$.

We claim that $$p_{A_M}( \Theta_M)\bullet p_{A_M}(\hat{\alpha})
\neq 0.$$ In expanding $p_{A_M}( \Theta_M)\bullet
p_{A_M}(\hat{\alpha})$ we obtain four terms, all in $\H_0(M)$ and
each one is represented by a single loop and a $\pm$ sign. These
four loops have free homotopy type $[\alpha.h],[h],[\alpha]$ and
$[1]$. Since there are at least three different classes, any
linear combination of them with coefficient $\pm1$ is nonzero.
Therefore $M$ has nontrivial extended loop products.
\end{proof}
\section{Seifert Manifolds with Boundary}\label{sec-Seifert-with-boundary}
By torus decomposition, Seifert manifolds with
incompressible\footnote{$\pi_1$-injective.} boundary are among the
building blocks of $3$-manifolds. In this section we provide
examples of nontrivial extended loop product for various Seifert
manifold with boundary. The examples are organized according to
the number of boundary components and singular fibers.

Throughout this section $M$ is a compact oriented Seifert manifold
with orientable base surface unless otherwise it is stated (see
Corollary \ref{2-b-nonorientable}).\\
 $M$ has $b\geq 1 $ boundary components and
$p$ singular fibers and $p'$ denotes the number of the singular
fibers of multiplicity
 greater than $2$. We assume  \footnote{This assumption will be justified later, in
section 10, when we recall the torus decomposition.} that
$$b+p\geq 3.$$
$g$ denotes the genus of the base surface of $M$. In the proofs of
the statements in this section $g$ plays no role. So we assume
that
$$g=0.$$

We present two main ideas in constructing nontrivial extended loop
products for Seifert manifolds with boundary,
\begin{center}
 \emph{Two curves argument} and \emph{Chas' figure eight argument}.
\end{center}
The following table indicates the cases where we apply these
arguments.
\begin{center}
\begin{tabular}{|c|c|}
  \hline
     & argument  \\
  \hline
$p+b\geq 4$& \small Two curves argument   \\
 \hline
$p'+b\geq 3$ &\small Chas' figure eight argument \\
  \hline
\end{tabular}
\end{center}

\bigskip

\subsection{$\mathbf{p+b \geq3}$ (Two curves argument)}\label{large-seifert}

\begin{proposition}\label{more-3-fiber}
Let $M$ be a compact oriented Seifert manifold with $p>2$ and
$b\geq 1$. Then $M$ has nontrivial extended loop products.
\end{proposition}
\begin{proof}
Let $c_1$, $c_2$ and $c_3$ be the generators of $\pi_1(M)$
corresponding to three singular fibers. Consider a simple curve
$\gamma_1$, on the base surface, away from singular fibers and
representing the free homotopy class $[c_1c_2]$. Similarly,
consider a simple curve $\gamma_2$ on the base surface, away from
the singular fibers and
 representing the free homotopy class $[c_2c_3]$. Moreover,
$\gamma_2$ can be chosen such that it has exactly $2$ intersection
points with $\gamma_1$.

Since $\gamma_1$ is away from the singular point, the fibration
can be trivialized over $\gamma_1$ (see Figure
\ref{3-sing-fibre}). Therefore, we obtain a map
$$
f:\T ^2 \rightarrow M
$$
from the torus to $M$, where $f(0,\cdot)=\gamma_1$.

$f$ gives rise to a homology class $\Gamma_1 \in \H _1(M)$, by
declaring $f(t,0)$, $t \in \S ^1$, as the set of marked points.
The loop passing through $f(t,0)$ is $f(t, s)$, $s\in S^1$ and has
the free homotopy type $[c_1c_2]$. Also, the loop $\gamma_2$ gives
rise to a homology class $\Gamma_2 \in \H _0 (M)$.
\begin{center}
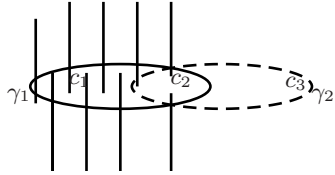
\begin{figure}[h]
\begin{pspicture}(1,0)(15,5)
%\psgrid(0,0)(15,5)
\psellipse[linewidth=1.0pt](6.5,2.5)(2.7,0.7)
\psellipse[linewidth=1.0pt,linestyle=dashed](9.5,2.5)(2.7,0.7)
%\psellipse[linewidth=1.0pt,linecolor=green](6.9,2.5)(1.2,3.3)
\psline[linewidth=1pt](8,2.8)(8,5)
\psline[linewidth=1pt](7,2.5)(7,5)
\psline[linewidth=1pt](6,2.3)(6,5)
\psline[linewidth=1pt](5,2.3)(5,5)
\psline[linewidth=1pt](4,2)(4,4.5)
\psline[linewidth=1pt](4.5,2.4)(4.5,0)
\psline[linewidth=1pt](4.5,2.9)(4.5,0)
\psline[linewidth=1pt](8,2.3)(8,0)
\psline[linewidth=1pt](6.5,2.9)(6.5,0)
\psline[linewidth=1pt](5.5,2.9)(5.5,0)
\psline[linewidth=1pt](4.5,2.9)(4.5,0)

%\pscurve[linewidth=1.0pt](8.5,2)(6.5,2)(7,3)(9,2.5)(8.5,3)
%\pscurve[linewidth=1.0pt,linestyle=dashed](9.95,2.4)(8.5,2)(8,2.5)(9,3)(11,2)(11.5,2.5)(11,3)(10.2,2.7)
%\pscurve[linewidth=1.0pt](6.9,4)(7.2,2.5)(6.8,1.1))(0.65,0.2)
\rput[b]{0}(5.3,2.5){\small $c_1$}\rput[b]{0}(3.5,2){\small
$\gamma_1$}\rput[b]{0}(12.5,2){\small $\gamma_2$}
\rput[b]{0}(8.3,2.5){\small $c_2$} \rput[b]{0}(11.7,2.4){\small
$c_3$}
%\rput[b]{0}(3,1){\small Set of marked points}
%\rput[b]{0}(3,0.7){of $H$}
%\psline{->}(3,1.4)(5.5,2.4)
\end{pspicture}
\caption{$\gamma_1$ and $\gamma_2$ and a trivialized fibration
over $\gamma_1$}\label{3-sing-fibre}
\end{figure}
\end{center}

 We claim that

$$
p_{A_M}(\Delta \Gamma_1)\bullet p_{A_M}(\Delta \Gamma_2) \neq 0
\in \H _0(M)
$$
which proves the proposition.

The sets of marked points of $\Delta \Gamma_1\in \H_2(M)$ and
$\Delta \Gamma_2 \in \H_1 (M)$ intersect at two points. Each point
contributes four terms in the expansion of $p_{A_M}(\Delta
\Gamma_1)\bullet p_{A_M}(\Delta \Gamma_2)$. All in all eigth terms
emerge while calculating $p_{A_M}(\Delta \Gamma_1)\bullet
p_{A_M}(\Delta \Gamma_2)$. Each is represented by a conjugacy
class in $\hat{\pi}_1(M)$ and a $\pm $ sign. These conjugacy
classes are

\begin{center}
$ [c_1c_2c_3c_2], [c_1c_2],[c_2c_3], [1] $\\ and\\
$ [c_1c_2^2c_3],[c_1c_2],[c_2c_3],[1]. $
\end{center}
\begin{center}
\begin{figure}[h]
\begin{pspicture}(-2,0)(15,10)
%\psgrid(0,0)(15,10)
\psellipse[linewidth=1.0pt](1.5,2.7)(0.7,2.5)
\psellipse[linewidth=1.0pt,linestyle=dashed](1.5,5.2)(0.7,2.5)
\rput[b]{0}(1.5,1.0){\small $c_1$} \rput[b]{0}(1.5,3.8){\small
$c_2$} \rput[b]{0}(1.5,6.5){\small $c_3$}
\rput[b]{0}(0.86,3.8){\small $\bullet$}
\rput[b]{0}(2.1,3.8){\small $\bullet$} \psline{->}(1,4)(4.8,6.8)
\pscurve[linewidth=1.0pt](5,7)(5.5,6.5)(5.8,7)
\pscurve[linewidth=1.0pt,linestyle=dashed](5.8,7)(5.5,7.5)(5,7)(5,5.5)(5.5,4.8)(6.2,5.5)(6.5,7)
\pscurve[linewidth=1.0pt](6.5,7)(6.2,8.5)(5.5,9)(5,7)
\rput[b]{0}(5.5,5.5){\small $c_1$} \rput[b]{0}(5.5,7){\small
$c_2$} \rput[b]{0}(7.5,7){\small $[c_1c_2c_3c_2]$}
\rput[b]{0}(5.5,8){\small $c_3$}
\pscurve[linewidth=1.0pt](5.4,2)(5.7,1.5)(6.2,2)
\pscurve[linewidth=1.0pt,linestyle=dashed](4.7,2)(5,0.5)(5.5,0)(6,0.5)(6.2,2)(5.7
,2.5)(5.4,2)
\pscurve[linewidth=1.0pt](6.2,2)(6,3.7)(5.5,4.2)(4.7,2)
\rput[b]{0}(5.7,0.5){\small $c_1$} \rput[b]{0}(5.7,2){\small
$c_2$} \rput[b]{0}(7.5,2){\small $[c_1c_2^2c_3]$}
\rput[b]{0}(5.7,3){\small $c_3$} \psline{->}(2.2,3.8)(4.6,2)
\end{pspicture}
\caption{Two curves argument}\label{2-curve-arg}
\end{figure}
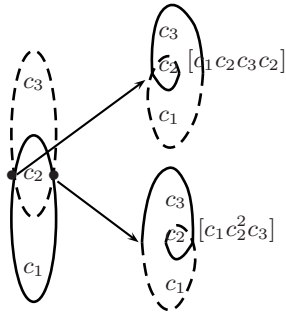
\end{center}
Thus in calculating $p_{A_M}(\Delta \Gamma_1)\bullet
p_{A_M}(\Delta \Gamma_2)$ $\mod 2$, only two terms remain, namely
(see Figure \ref{2-curve-arg})

$$[c_1c_2^2c_3] \quad \text{and} \quad [c_1c_2c_3c_2].$$

To prove the claim it is sufficient to show that these two
conjugacy classes are different. For this, consider the group
$$H=\langle c_1,c_2,c_3,d| c_1c_2c_3d_1=1, c_i^{\alpha_i}=1, 1\leq i \leq 3 \rangle,$$
where $d_1$ is a generator corresponding to a boundary component.
Indeed $H$ is the free product of the groups
$\Z_{\alpha_1}*\Z_{\alpha_2}*\Z_{\alpha_3}$.

Consider the homomorphism $\phi: \pi_1(M)\rightarrow H$ which is
 identity on $c_1,c_2,c_3,d_1$ and sends the other generators to the trivial element of $H$.
By the presentation (\ref{seifert-oriented}) $\phi$ is
well-defined.

Note that the image  of the conjugacy classes $[c_1c_2^2c_3]$ and
$[c_1c_2c_3c_2]$ under the $\phi $ are different in $\hat{H}$  as
$c_1c_2^2c_3$ and $c_1c_2c_3c_2$ are cyclically different reduced
words. Therefore $[c_1c_2^2c_3]$ and $[c_1c_2c_3c_2]$ are
different conjugacy classes in of $\pi_1(M)$.

\end{proof}

\begin{remark}
\label{remark-boundary} \emph{Indeed, in the proof of the
proposition above one can replace a singular fiber by a boundary
component or in other words a boundary component works as well as
a singular fiber of multiplicity zero.
 To be more explicit, the curves $\gamma_1$ and $\gamma_2$
 can go around a boundary component of the base surface instead of the singular fiber.
 In the free product $H$, one of finite cyclic groups is replaced by a infinite cyclic group.
  The rest of the proof remains the same.}
\end{remark}

\begin{corollary}\label{coro-3-fiber}
Let $M$ be a compact oriented Seifert manifold with $p$ singular
fibers and $q$ boundary components. If $p+b>3$ then $M$ has
nontrivial extended loop products.
\end{corollary}
\begin{proof}
 If $p\geq 3$ this is just the statement of Proposition \ref{more-3-fiber}.
  In the case $p< 3$, by Remark \ref{remark-boundary},
  in the proof of Proposition \ref{more-3-fiber} one replaces the missing singular fiber with one of
 the extra boundary components.
\end{proof}
\subsection{$p'+b\geq3$, (Chas' figure eight argument)}
\label{small-fiber}

 Now we turn to compact oriented Seifert manifolds with the total
number of singular fibers and boundary components equals $3$.

The idea of the proof of this case is due to Moira Chas. She found
this idea while reformulating a conjecture of Turaev (see
\cite{Chas}) on Lie bi-algebras of surfaces, characterizing non
self-intersecting closed curves.

\begin{proposition}
\label{2s-1b} Let $M$ be an oriented Seifert manifold with $p=2$
singular fibers and $b=1$ boundary component. Suppose that none of
its singular fibers has multiplicity $2$. Then $M$ has nontrivial
extended loop products.
\end{proposition}
\begin{proof}
Let $c_1$ and $c_2$ be generators of $\pi_1(M)$ corresponding to
the singular fibers. Let $\gamma_1$ be a smooth curve on the base
surface, with one self intersection point, away from singular
fibers and representing the free homotopy class $[c_1c_2^{-1}]$.
Similarly, consider a curve $\gamma_2$ on the base surface, away
from the singular fibers, with one intersection point and
representing the free homotopy class $[c_1^{-1}c_2]$. Moreover,
$\gamma_2$ can be chosen such that it has exactly $2$ intersection
points with $\gamma_1$.

The fibration can be trivialized over the $\gamma_1$ since it is
away from the singular point. Therefore, we obtain a map
$$
f:\T ^2 \rightarrow M
$$
from torus to $M$, where $f(\cdot,0)=\gamma_1$. This map gives
rise to a homology class $\Gamma_1 \in \H _1(M)$, by declaring
$f(t,0)$, $t \in \S ^1$, as the set of marked points. The loop
passing through $f(t,0)$ is $f(t,s)$, $s\in S^1$. All the loops of
this family have the free homotopy type $[c_1c_2^{-1}]$. The loop
$\gamma_2$ gives rise to a homology class $\Gamma_2 \in \H _0
(M)$. We claim that

$$
p_{A_M}(\Delta \Gamma_1)\bullet p_{A_M}(\Delta \Gamma_2) \neq 0
\in \H _0(M)
$$
which proves the proposition.

The set of marked points  of $\Delta \Gamma_1\in \H_2(M)$ and
$\Delta \Gamma_2 \in \H_1 (M)$ intersect exactly at two points.
Each point contributes four terms in the expansion of
$p_{A_M}(\Delta \Gamma_1)\bullet p_{A_M}(\Delta \Gamma_2)$. All in
all $8$ terms emerge in the expansion of $p_{A_M}(\Delta
\Gamma_1)\bullet p_{A_M}(\Delta \Gamma_2)$. Each ones of them is
represented by a conjugacy class in $\hat{\pi}_1(M)$ and a $\pm $
sign. These conjugacy classes are
\begin{center}
$ [c_1c_2c_1^{-1}c_2^{-1}], [c_1^{-1}c_2],[c_1c_2^{-1}], [1] $\\
and\\
$ [c_1c_2^{-1}c_1^{-1}c_2], [c_1^{-1}c_2],[c_1c_2^{-1}], [1]. $
\end{center}

\begin{center}
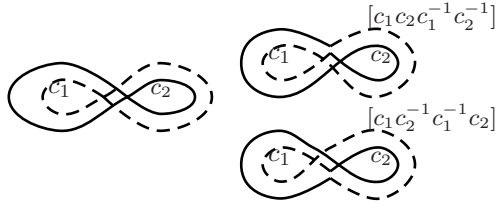
\begin{figure}[h]
\begin{pspicture}(3,0)(15,5)

%\psgrid(0,0)(15,5)
\pscurve[linewidth=1.0pt](6.5,3)(4.5,4)(3,3)(4.5,2)(6.5,3)(7.5,3.5)(8.5,3)(7.5,2.5)(6.5,3)
\pscurve[linewidth=1.0pt,linestyle=dashed](5.8,3)(4.5,3.5)(4,3)(4.5,2.5)(5.8,3)(7.5,4)(9,3)(7.5,2)(5.8,3)
\pscurve[linewidth=1.0pt](12.5,4.5)(11,5)(10,4.5)(10,3.5)(11,3)(14,4.5)(14.5,4)(14,3.5)(12.5,4.3)
\pscurve[linewidth=1.0pt,linestyle=dashed](12.5,4.3)(11,3.5)(10.5,4)(11,4.5)(12,4.3)(13.8,3)(15,4)(14,5)(12.5,4.5)
\pscurve[linewidth=1.0pt](12.5,0.8)(14,1.5)(14.5,1)(14,0.5)(12,1.5)(11,2)(10,1.5)(10,0.5)(11,0)(12.5,0.6)
\pscurve[linewidth=1.0pt,linestyle=dashed](12.5,0.6)(14,0)(15,1)(14,2)(12.5,1.5)(11,0.5)(10.5,1)(11,1.5)(12.5,0.8)
\rput[b]{0}(4.5,3){$c_1$}\rput[b]{0}(11,1){$c_1$}\rput[b]{0}(11,4){$c_1$}
\rput[b]{0}(7.5,3){\small $c_2$}\rput[b]{0}(14,1){\small
$c_2$}\rput[b]{0}(14,4){\small $c_2$} \rput[b]{0}(15.5,5){\small
$[c_1c_2c_1^{-1}c_2^{-1}]$}\rput[b]{0}(15.5,2){\small
$[c_1c_2^{-1}c_1^{-1}c_2$]}
\end{pspicture}
\caption{Figure eight argument}\label{fig-gamma}
\end{figure}
\end{center}

Calculating $p_{A_M}(\Delta \Gamma_1)\bullet p_{A_M}(\Delta
\Gamma_2)\mod 2$, only two terms remain (see Figure
\ref{fig-gamma}), namely

\[
[c_1c_2c_1^{-1}c_2^{-1}]\quad  \text{and} \quad
[c_1c_2^{-1}c_1^{-1}c_2].
\]

To prove the claim it is sufficient to show that these two
homology classes are different. For this, consider
 the subgroup
 $$H=\langle c_1,c_2,d_1| c_1c_2d_1=1 ,c_i^{\alpha_i}=1, i=1,2      \rangle$$
where $d_1$ corresponds to one of the boundary components in
presentation (\ref{seifert-oriented}). Indeed $H$ is the free
product of cyclic groups generated by $c_1$, $c_2$ and $d_1$.

Consider the homomorphism $\phi: \pi_1(M)\rightarrow H$ which is
identity on $c_1,c_2,d_1$ and sends the other generators of
$\pi_1(M)$ in presentation (\ref{seifert-oriented}) to the
identity element of $H$. It follows from presentation
(\ref{seifert-oriented}) that $\phi$ is well defined.

Since $\alpha_i\neq 2$, $i=1,2$, then $c_i\neq c_i^{-1} \in H$.
Thus the images of $ c_1c_2c_1^{-1}c_2^{-1} $ and
$c_1c_2^{-1}c_1^{-1}c_2$ under $\phi$ in the free product $H$ are
not conjugate as they are cyclically different. Therefore the
conjugacy classes $[c_1c_2c_1^{-1}c_2^{-1}]$ and
$[c_1c_2^{-1}c_1^{-1}c_2]$ in $\hat{\pi}_1(M)$ are different as
well.
\end{proof}
\begin{remark}\emph{
The idea of the proof above works unless there is a singular fiber
of multiplicity $2$ and this case needs a special
treatment}.\end{remark}
\begin{corollary}
\label{3boundary} Let $M$ be a compact oriented Seifert manifold
with $p'+b\geq 3$. Then $M$ has nontrivial extended loop products.
\end{corollary}

\begin{proof}
If $p'\geq 2$ then it is just Proposition \ref{2s-1b}. Otherwise,
by Remark \ref{remark-boundary} a boundary component works as well
as a singular fiber with multiplicity zero. The proof is the same
as of the one of Proposition \ref{2s-1b} except that one has to
replace one of the generators of $\pi_1(M)$ contributed by a
singular fiber to a generator corresponding to a boundary
component.
\end{proof}
\begin{corollary}
\label{2-b-nonorientable}
 Let $M$ be a compact oriented Seifert manifold with a non oriented
 base and $b\geq 2$ boundary components, then $M$ has nontrivial
extended loop products.
 \end{corollary}
\begin{proof}
We can modify the proof of Proposition \ref{2s-1b} considering
that by Remark \ref{remark-boundary} a boundary component serves
our purposes as well as a singular fiber.  Note that the fibration
is oriented once restricted  to the boundary components. So it can
be trivialized over a simple curve representing the free homotopy
type $[d_1d_2^{-1}]$ where $d_1$ and $d_2$ are generators of
$\pi_1(M)$ contributed by the two boundary components (see the
presentation \ref{seifert-nonoriented}). The rest of  the proof is
similar to the one of Proposition \ref{2s-1b}. Calculating mod 2
we get two terms

$$
[d_1d_2d_1^{-1}d_2^{-1}] \quad \& \quad [d_1d_2^{-1}d_1^{-1}d_2].
$$

 If $a_1$ is a
generators of $\pi_1(M)$ contributed by a cross cap (see the
presentation \ref{seifert-nonoriented}) then consider the group
\begin{equation*}
H=\langle d_1,d_2,a_1| d_1d_2a_1^2=1\rangle\simeq  \langle d_1
\rangle * \langle a_1 \rangle
\end{equation*}
and the homomorphism $\varphi:\pi_1(M)\rightarrow H$ which is
identity on $d_1,d_2,a_1$ and sends the other generators of
$\pi_1(M)$. The image of the two conjugacy classes above in the
free product $H$ is $[d_1a_1^2d_1^{-1}a_1^{-2}] \quad \& \quad
[d_1a_1^{-2}d_1^{-1}a_1^2]$ which are different. Therefore the
conjugacy classes above are distinct conjugacy classes of
$\pi_1(M)$.

\end{proof}

\section{Connected Sum of $3$-Manifolds}

In this section we consider the connected sum of $3$-manifolds.
 We shall show how one can obtain nontrivial loop products in this kind of
 manifold. For this part, the author has benefited from conversations with
 numerous colleagues \cite{Bend}.

\begin{proposition}
\label{connected-sum}
 If $M_1$ and $M_2$ are two 3-manifolds such that
$\pi_1(M_i)\neq {1}$, $i=1,2$, then $M_1\#M_2$ has nontrivial
extended loop products.
\end{proposition}
\begin{proof}
Let $m_i\in M_i$, $i=1,2$ be two points which are identified in
the connected sum $M_1\#M_2$ and considered as the base point $m$
of $M_1$, $M_2$ and $M_1\#M_2$. Let $g_1\in \pi_1(M_1)$ and $g_2
\in \pi_1(M_2)$ be two non-trivial elements. Consider the loops
$\gamma_i$ in $M_i$ with marked point $m$ and representing $g_i$,
$i=1,2$.

 By composing these two loops at $m$, we obtain
a loop $\gamma_1\gamma_2$ representing $g_1g_2 \in \pi_1(M_1\#
M_2)\simeq \pi_1(M_1)*\pi_1(M_2)$ and it has the free homotopy
type $[g_1g_2]\in \hat{\pi}_1(M)$ (see Figure \ref{conn-sum}).
Notice that $\gamma_1\gamma_2$ also represents an element of
$\mathbb{H}_0(M_1\#M_2)$. We denote this homology class by
$g_1\#g_2$.

Consider the embedding $\beta:S^2 \hookrightarrow M_1\#M_2$ where
the connected sum happens (see Figure \ref{conn-sum}). This gives
rise to an element in $\pi_2(M_1\#M_2) \simeq \pi_1(\Omega(
M_1\#M_2))$ and therefore an element in
$\pi_1(\mathbb{L}(M_1\#M_2))$. The image of this element under the
Hurewicz map is an element of $\mathbb{H}_1(M_1\#M_2)$. By
composing the loops of this $1$-cycle of loops with a fixed loop
whose marked point is $m$ and represents a nontrivial element $h$
in $\pi_1(M_1\#M_2)$, we can obtain a new element in
$\mathbb{H}_1(M_1\#M_2)$. We denote this element by
$\beta_{M_1\#M_2}$ and the homotopy type of its loops is $[h]\in
\hat{\pi}_1(M_1\#M_2)$.

\begin{center}
\begin{figure}[h]
\begin{pspicture}(0,-1)(15,5)
 \pscircle[linewidth=1.0pt](7,2.5){2.5}
\psellipse[linewidth=1.0pt](7,2.5)(0.4,2.5)
\rput[b]{26.51}(0.8,-3.4){\psellipse[linewidth=1.0pt](7,2.5)(0.3,2.23)}
\rput[b]{-26.51}(0.6,2.9){\psellipse[linewidth=1.0pt](7,2.5)(0.3,2.23)}
\rput[b]{11.76}(0.06,-1.49){\psellipse[linewidth=1.0pt](7,2.5)(0.35,2.45)}
\rput[b]{-11.76}(0.06,1.44){\psellipse[linewidth=1.0pt](7,2.5)(0.35,2.45)}
\rput[b]{38.659}(1.9,-4.8){\psellipse[linewidth=1.0pt](7,2.5)(0.3,1.920)}
\rput[b]{-38.659}(1.2,3.9){\psellipse[linewidth=1.0pt](7,2.5)(0.3,1.920)}
\psccurve[linewidth=1.0pt,linestyle=dashed](3,3)(5.8,4)(10,3.8)(11,3)(5.8,2)(3,3)
\pscurve[linewidth=1.0pt](2,5.3)(7,5)(9,5.3)(11,4.8)(12.5,5.3)
\pscurve[linewidth=1.0pt](2,-1.5)(5,-1.2)(7,0)(9,-1.5)(11,-1)(12.5,-1.3)
\pscurve[linewidth=1.0pt,linestyle=dashed](7,0)(5,1.2)(3,1)(3,0.8)(7,0)
\pscurve[linewidth=1.0pt,linestyle=dashed](7,0)(9,1.2)(11,1)(11,0.8)(7,0)
\rput[b]{0}(5.8,3.9){\small$\bullet$}
\rput[b]{0}(5.8,1.9){\small$\bullet$}
\rput[b]{0}(2,2.3){\Large$M_1$} \rput[b]{0}(13,2.5){\Large$M_2$}
\rput[b]{0}(2.7,1){\Large$\gamma_1$}
\rput[b]{0}(11.5,1){\Large$\gamma_2 $}
\rput[b]{0}(4,1.8){\Large$\gamma$}
\rput[b]{0}(12.5,4){\Large$g_1\#g_2$}
\psline{->}(11.7,4.5)(9.5,4.2)
\end{pspicture}
\caption{Loop product in connected sum}\label{conn-sum}
\end{figure}
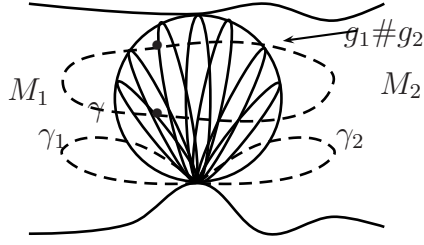
\end{center}

\noindent We claim that
$$p_{A_M}(\Delta(g_1\#g_1))\bullet p_{A_M}(\Delta(\beta_{M_1\#M_2})) \neq 0$$
which proves that $M$ has nontrivial extended loop product. To
compute
\begin{equation}
\begin{split}
p_{A_M}(\Delta(g_1\#g_2))\bullet
p_{A_M}(\Delta(\beta_{M_1\#M_2}))=(\Delta(g_1\#g_2)-p_M(\Delta(g_1\#g_2)))\\
                             \bullet(\Delta(\beta_{M_1\#M_2})-p_M(\Delta(\beta_{M_1\#M_2})))
\end{split}
\end{equation}
 one has to homotop
$\gamma_1\gamma_2$ to a simple curve $\gamma$ which intersects
 the embedded $S^2$ in $M_1\#M_2$  at two different points. So
in calculating the expansion above we get eight terms, four terms
for each intersection point of $\gamma$ and the embedded
$2$-sphere. Each one of these eight terms, which are elements of
$\H _0(M)$, is represented by a single loop and a $\pm$ sign.
These loops have free homotopy types,
\begin{equation*}
[g_1g_2h],[g_1g_2],[h],[1],[g_2g_1h],[g_2g_1],[h],[1].
\end{equation*}

So calculating $p_{A_M}(\Delta(g_1\#g_2))\bullet
p_{A_M}(\Delta(\beta_{M_1\#M_2}))$ mod  2, only two terms remain,
$$[g_1g_2h]\& [g_2g_1h].$$

\noindent If we choose $h=g_1g_2$ in $\pi_1(M_1\#M_2)$, the two
conjugacy classes
$$[g_1g_2h]=[g_1g_2g_1g_2] \text{ and } [g_2g_1h]=[g_1g_1^2g_2]=[g_2^2g_1^2]$$ are different as two
cyclically reduced sequences $(g_1,g_2,g_1,g_2)$ and
$(g_2^2,g_1^2)$ are not cyclic permutation of one another.

 This proves that
for $h=g_1g_2$, $ p_{A_M}(\Delta(g_1\#g_2))\bullet
p_{A_M}(\Delta(\beta_{M_1\#M_2}))$ is nonzero$\mod 2$, therefore
it is nonzero in $\H_0(M_1\#M_2)$.
\end{proof}

\section{An Injectivity Lemma}
Suppose that $M$ is a closed oriented manifold and $\mathcal{T}$
is a collection of incompressible \footnote{$\pi_1$-injective}
separating tori in $M$. Let $M\setminus \mathcal{T}=M_1\cup
M_2...\cup M_n$. We associate with $(M,\mathcal{T})$ a tree of
groups $(G,T)$ (See Appendix A for the definitions).

The vertices of $T$ are in one-to-one correspondence with $M_i$'s. Two vertices $v_i$
and $v_j$ are connected by an edge if the corresponding $M_i$ and $M_j$ are glued along a torus.
With each vertex $v_i$ we associate the fundamental group of $\pi_1(M_i)$ and with each edge
we associate the fundamental group of the corresponding torus. As the tori are incompressible
there are two monomorphisms from each edge group to the vertex groups of the corresponding
vertices. Namely if $e$ is an edge (a torus in $M$) and $v$ a vertex of $e$
(some $M_i$ s.t. $T\subset \partial \overline{M_i}$ ) then there is an homomorphism
$$
f_e^v:G_e\rightarrow G_v
$$
which is indeed the map induced by inclusion
$$
i:\pi_1(T) \rightarrow \pi_1(\overline{M_i})\simeq \pi_1(M_i).
$$

By Van-Kampen theorem, $\pi_1(M)$ is the amalgamation of $G_v$'s along $G_e$'s (see Appendix II, tree of groups).

Combining all these with Lemma \ref{lemma-tree} in Appendix A, we
have,
\begin{lemma}(Injectivity Lemma)
\label{injection-lem} Let $M$, $\mathcal{T}$ and $M_i$'s be as
above. Suppose that $[a]$ and $[b]$ are two distinct conjugacy
classes in $\pi_1(M_i)$ such that $a$ is not conjugate to the
image of any element of the edge groups in $\pi_1(M_i)$. Then $a$
and $b$ represent distinct conjugacy classes of $\pi_1(M)$.
\end{lemma}

The injectivity lemma will be very useful to show that certain
examples of nontrivial loop products in a submanifold of $M$ give
rise to nontrivial loop product in $M$.

\section{Gluing Seifert Manifolds along Tori}

In this section we consider the gluing of the Seifert manifolds
along tori. We prove that the result of gluing Seifert manifolds
along incompressible \emph{separating} tori has nontrivial
extended loop products except certain classes for which there are
double covers that have nontrivial extended loop products.

The idea of the proof is to apply either the \textit{two curves
argument} or the \textit{figure eight argument} to a Seifert piece
of the gluing that fulfills the condition of the arguments. Then
we show that the examples of nontrivial loop product in the piece
give rise to a nontrivial example of loop product in the manifold,
for that we need the injectivity lemma, Lemma \ref{injection-lem}.

When we cannot apply either the two curves argument or the figure
eight argument we consider an appropriate double cover of the
manifold which either has a non separating torus or one of the two
arguments can be applied to one of the pieces in its torus
decomposition
 \footnote{See Section 10 for the definition of torus decomposition}.

Before stating the result and demonstrating the proof we introduce
the following notation.

\noindent \textbf{Convention and Notation:}
 Let $M$ be a closed 3-manifold and $\mathcal{T}\neq \emptyset$
be a collection of tori in  $M$ and
 $$M\setminus \mathcal{T}=M_1 \cup M_2 \cdots  \cup M_ n,$$
where all the $\overline{M_i}$'s are oriented Seifert manifolds
with incompressible boundary\footnote{We say that the torus
decomposition of $M$ has only Seifert pieces.}.

 Let $b_i \geq 1$ and $p_i$ denote respectively the number of boundary
components and singular fibers of $\overline{M_i}$ and $A_i\subset
\Z^+$ is the set of all the multiplicities of the singular fibers
of $\overline{M_i}$. The base surface of $M_i$ is denoted $S_i$
and the genus of $S_i$ is $g_i$ if $S_i$ is oriented otherwise
 $g_i$ is the number of cross caps in $S_i$.

We assume that all tori in $T$ are separating since in Section
\ref{nonseparating-section} we proved that 3-manifolds with non
separating torus have nontrivial extended loop products.

We assume\footnote{see section 10 for justification.} that for all
$i$
$$
b_i+p_i\geq 3
$$
if $S_i$ is orientable.
 Also we assume that $$g_i=0$$  if $S_i$ is orientable, as if $g_i>
 0$ then there is a non separating torus in $M$ and this case has already been studied
 in Section \ref{nonseparating-section}.
\begin{proposition}
\label{seifert-glued}
 Let $M$ be a closed 3-manifold  described as
above. Then $M$ or a double cover of $M$ has nontrivial extended
loop products.

\end{proposition}
\begin{proof} $M$ belongs to one of the following categories:
\begin{numlist}
\item $b_i\geq 3$ and $S_i$ is orientable for some $i$. In this
case we use the figure eight argument. Note that if $b_i >3$ we
can use the two curves argument. \item $b_i\geq 2$ and $S_i$ is
not orientable for some $i$, for this case we also use the figure
eight argument .
 \item $p_i+b_i
>3$ and $S_i$ orientable for some $i$, we can apply two curves argument.

\item $b_i+p_i=3$, $1\leq b_i\leq2$, $S_i$ orientable and $2\notin
A_i$ for some $i$. In this case we use figure eight argument.

\item For all $i$, $b_i+p_i=3$, $1\leq b_i\leq2$, $S_i$ orientable
and $2\in A_i$ for all $i$.
 For this case we consider certain double cover of $M$.

\item $b_i+p_i=3$ and $2\in A_i$ if $S_i$ is orientable otherwise
$b_i=1$ and $S_i$ is not orientable. In this case also we consider
a double cover.

\end{numlist}

\noindent We verify the statement for each cases.

\begin{numlist}

\item \textbf{$b_i \geq3$ and $S_i$ orientable for some $i$:}

Suppose that $b_1\geq 3$. By Corollary \ref{3boundary}
$\overline{M_1}$ has nontrivial extended loop products.
 Consider the homology classes $\Gamma_1$ and $\Gamma_2 \in \H_1(\overline{M_1})$
constructed in the proof of Corollary \ref{3boundary}. $\Gamma_1$
and $\Gamma_2$ can be regarded as the elements of  $\H_*(M)$. We
claim that
$$
p_{A_M}(\Delta \Gamma_1)\bullet p_{A_M}(\Delta \Gamma_2) \neq 0
\in \H _0(M).
$$

In calculating the  product above, four terms emerge, all of them
in $\H_0(M)$.
 To prove the product is nonzero it is sufficient to show it is nonzero mod 2. Having done the
 calculation mod 2, just like in the proof of Proposition \ref{2s-1b} and Corollary \ref{3boundary},
 only two terms survive, $$[d_1d_2d_1^{-1}d_2^{-1}] \text{ and } [d_1d_2^{-1}d_1^{-1}d_2],$$

\noindent where $d_1$ and $d_2$ are the generators contributed by
two boundary components (see presentation
(\ref{seifert-oriented})).

By the proof of Corollary \ref{3boundary} (indeed Proposition
\ref{2s-1b})
 we know these two are distinct conjugacy classes of $\pi_1(M_1)$. So if we show that
$d_1d_2d_1^{-1}d_2^{-1}$ is not conjugate to any element of the
fundamental group of one of the components of $\partial\overline{
M_1}$, then by injectivity lemma $[d_1d_2d_1^{-1}d_2^{-1}]$ and
$[d_1d_2^{-1}d_1^{-1}d_2]$ are different as the conjugacy classes
of $\pi_1(M)$.

Suppose that

\begin{equation}
\label{eq-3b-conju} [d_1d_2d_1^{-1}d_2^{-1}]=[h^rd^s]
\end{equation}

\noindent where $h$ is the generator of $\pi_1(M_1)$ that
corresponds to a normal fiber of $M_1$ and $d$ is a generator of
$\pi_1(M)$ contributed by a boundary component.

Consider the group $H$
$$
H=\langle d_1,d_2,d_3|d_1d_2d_3=1\rangle\simeq \langle d_1
\rangle*\langle d_2\rangle
$$
and the homomorphism
$$\phi:\pi_1(M_1)\rightarrow H$$
which is identity on $d_1,d_2, d_3$ and trivial
 on  all other generators of $\pi_1(M_1)$.
 It follows from the presentation (\ref{seifert-oriented}) of
$\pi_1(M_1)$ that $\phi$ is well-defined.

 After applying $\phi$ to (\ref{eq-3b-conju}) we get

$$
[d_1d_2d_1^{-1}d_2^{-1}]=[1]
$$
if $d\neq d_1,d_2,d_3$ otherwise

$$
[d_1d_2d_1^{-1}d_2^{-1}]=[d^s].
$$

The first case is impossible as $d_1d_2d_1^{-1}d_2^{-1}$ is a
cyclically reduced word in the free product $\langle d_1 \rangle*
\langle d_2\rangle$ thus it represents a nontrivial conjugacy
class.

In the latter case, If $d=d_1$ (or $d_2$) then

$$
[d_1d_2d_1^{-1}d_2^{-1}]=[d_1^s],
$$
 and this is not possible either as $d_1d_2d_1^{-1}d_2^{-1}$ is
 cyclically reduced of length 4 and $d_1^s$ (or $d_2^s$) has
 length 1, hence they represent different conjugacy classes.

 If $d=d_3=d_2^{-1}d^{-1}$ then we have
$$
[d_1d_2d_1^{-1}d_2^{-1}]=[(d_1d_2)^{-s}].
$$
Considering the length of the word we have to have $s=\pm 2$. In
either case this inequality does not hold as $d_1\neq d_1^{-1}$
and $d_2\neq d_2^{-1}$($H$ is a free group). Therefore these two
conjugacy classes are always distinct. This proves the claim.

\item \textbf{$b_i \geq 2$ and $S_i$ is not orientable:}\\
 Similar to the previous case we use the figure eight argument in $M_i$.
 Suppose that $b_1\geq 2$ and $S_1$ is not orientable. By
Corollary \ref{2-b-nonorientable} $\overline{M_1}$ has nontrivial
extended loop products. The homology classes $\Gamma_1$ and
$\Gamma_2$ in $\H_*(\overline{M_1})$ constructed in the proof of
Corollary \ref{2-b-nonorientable} can be regarded as homology
classes of $\H_*(\overline{M})$. We claim that
$$
p_{A_M}(\Delta \Gamma_1)\bullet p_{A_M}(\Delta \Gamma_2) \neq 0
\in \H _0(M).
$$
Calculating$\mod 2$ we get two terms represented by the conjugacy
classes $[d_1a_1^2d_1^{-1}a_1^{-2}]$ and
$[d_1a_1^{-2}d_1^{-1}a_1^{2}]$ (see the proof of Corollary
\ref{2-b-nonorientable}). We must prove that these two conjugacy
are distinct as conjugacy classes of $M$. For that, by injectivity
lemma, it is sufficient to prove that one of them has no boundary
representative. We show that $[d_1a_1^2d_1^{-1}a_1^{-2}]$ has no
boundary representative. Suppose we have
\begin{equation}
\label{2b-nonori}
 [d_1a_1^2d_1^{-1}a_1^{-2}]=[h^rd^s].
\end{equation}
where $d$ is a generator of $\pi_1(M_1)$ contributed by a boundary
component and $h$ corresponds to the normal fiber. Consider the
group
$$
H=\langle d_1,d_2,a_1| d_1d_2a_1^{2}=1 \rangle\simeq \langle
d_1\rangle*\langle a_1\rangle \simeq \langle d_2\rangle *\langle
a_1 \rangle
$$
and the homomorphism $\varphi:\pi_1(M_1)\rightarrow H$ which is
identity on $d_1,d_2,a_1$ and trivial on the other generators of
$\pi_1(M_1)$. It follows from the presentation
(\ref{seifert-nonoriented}) that $\varphi$ is well defined.
Applying $\varphi$ to (\ref{2b-nonori}) we have
\begin{equation*}
[d_1a_1^2d_1^{-1}a_1^{-2}]=[1]
\end{equation*}
if $d\neq d_1,d_2$, otherwise
$$
[d_1a_1^2d_1^{-1}a_1^{-2}]=[d^s]
$$
in the free product $H\simeq\langle d_1\rangle*\langle
a_1\rangle$. The first case is impossible as
$d_1a_1^2d_1^{-1}a_1^{-2}$ is a cyclically reduced word of length
4 so it represents a nontrivial conjugacy class. If $d=d_1$ then
$$
[d_1a_1^2d_1^{-1}a_1^{-2}]=[d_1^s]
$$
which is again impossible in the free product $\langle d_1
\rangle*\langle a_1 \rangle$ as the cyclically reduced words
$d_1a_1^2d_1^{-1}a_1^{-2}$  and $d_1^s$ are different in the free
product $\langle d_1 \rangle* \langle a_1 \rangle$. If $d=d_2$
then
$$
[d_1a_1^2d_1^{-1}a_1^{-2}]=[(a_1^2d_1)^{-s}].
$$

Again similar to previous cases we must have $s=\pm 2$. Even this
is not enough as $d_1\neq d_1^{-1}$ and $a_1^2\neq  a_1^{-2}$.

 \item \textbf{$p_i+b_i \geq 4$ and $S_i$ orientable for some $i$
:}\\
Since case 1 deals with the case $b_i\geq 3$ for some $i$ so we
may assume that $b_i\leq 2$. So we are dealing with one of the
following cases:
\begin{romlist}
\item   $p_i>2$ for some $i$ \item $b_i=p_i=2$ for some $i$.
\end{romlist}

The proof of both cases are similar and use the two curves
argument. We just have to remember that a boundary component works
as a singular fiber of multiplicity $0$. Here we only present the
proof of  case $(i)$.

Suppose that $p_1>2$. Consider the homology classes $\Gamma_1$ and
$\Gamma_2$ in $\H_*(\overline{M_1})$
 introduced in the proof of Proposition \ref{more-3-fiber}. They can be considered as the elements of
 $\H_*(M)$. We claim that
$$
p_{A_M}(\Delta \Gamma_1)\bullet p_{A_M}(\Delta \Gamma_2) \neq 0
\in \H _0(M)
$$
The calculation is the same and the results is a sum of conjugacy
classes of $\pi_1(M)$. Computing mod 2 only two terms
$$[c_1c_2c_3c_2] \text{ and } [c_1c_2^2c_3].$$

We must prove that they are different conjugacy classes of
$\pi_1(M)$. By the proof of Proposition \ref{more-3-fiber} we know
that they are distinct as the conjugacy classes of $\pi_1(M_1)$.
If we show that at least one of them has no representative
 in the fundamental group of any of the boundary components of $M_1$
 then by injectivity lemma they are different as conjugacy classes of
 $\pi_1(M)$. Suppose that
\begin{equation}
\label{eq-3p-1} [c_1c_2c_3c_2]=[h^rd^s]
\end{equation}
where $h$ stands for the normal fiber and $d$ is the generator
coming from a boundary component.

Consider the group,
$$
H=\langle c_1,c_2,c_3|c_1c_2c_3=1,c^{\alpha_i}=1, 1\leq i\leq
3\rangle,
$$
where $\alpha_i$ is the multiplicity of the singular fiber
corresponding to $c_i$, and the homomorphism
$\phi:\pi_1(M_1)\rightarrow H$ which is identity on $c_1,c_2,c_3$
and sends the other generators of $\pi_1(M_1)$ to the trivial
element of $H$. Applying $\phi$ to (\ref{eq-3p-1}) we have
$$
[c_1c_2c_3c_2]=[1]
$$
or
$$
c_1c_2c_3c_2=1
$$
in $H$. Since $c_1c_2c_3=1$ in $H$ we conclude that $c_2=1$ which
is a contradiction. \item \textbf{$b_i+p_i=3$, $1\leq b_i\leq 2$,
$S_i$ orientable and $2\notin A_i$ for some $i$:}\\
Here we shall use the figure eight argument. Suppose $i=1$. By
assumption, one of the following holds:
\begin{center}$( b_1=1
\text{ and }p_1=2)$ or $(b_1=2 \text{ and } p_1= 1.)$ \end{center}

The proof is similar for both cases, considering that a boundary
component behaves just like a singular fiber of multiplicity $0$.
So we only present the proof of the first case.

Suppose $p_1=2$ and $b_1=1$. Consider the homology classes
$\Gamma_1$ and $\Gamma_2 \in \H_*(\overline{M_1})$ as constructed
in the proof of the Proposition \ref{2s-1b}. They can be regarded
as homology classes in $\H_*(M)$. We prove that
$$
p_{A_M}(\Delta \Gamma_1)\bullet p_{A_M}(\Delta \Gamma_2) \neq 0
\in \H _0(M)
$$
which proves that $M$ has nontrivial extended products.

Having done the same mod 2 calculation as the proof of Proposition
\ref{2s-1b} we get two conjugacy classes
\[
[c_1c_2c_1^{-1}c_2^{-1}]\quad  \text{and} \quad
[c_1c_2^{-1}c_1^{-1}c_2].
\]
Just like the previous cases, in order to show that these two
conjugacy classes are distinct, we must prove that at least one of
them has no representative in the fundamental group of the
boundary component.

Suppose that
\begin{equation}\label{eq-1b-2s}
[c_1c_2c_1^{-1}c_2^{-1}]=[h^rd^s]
\end{equation}
where $h$ is the generator corresponding to the normal fiber and
$d$ is the generator corresponding to a boundary component.

Consider the group
$$
H=<c_1,c_2,d_1|c_1c_2d_1=1, c_1^{\alpha_1}=c_2^{\alpha_2}=1>
$$
and the homomorphism $\phi: \pi_1(M_1)\rightarrow H$ which is
identity on $c_1,c_2,d_1$ and sends other generators of $\pi_1(M)$
to the trivial element of $H$. It follows from presentation
(\ref{seifert-oriented}) that $\phi$ is well-defined.

By applying $\phi$ to (\ref{eq-1b-2s}) we have
$$
[c_1c_2c_1^{-1}c_2^{-1}]=[d_1^s].
$$
in $H$, or
$$
[c_1c_2c_1^{-1}c_2^{-1}]=[(c_1c_2)^s].
$$
Note that $H$ is the free product $\Z_{\alpha_1}*\Z_{\alpha_2}$
with $c_1$ and $c_2$ as the generators of the corresponding
factors.

So if $c_1c_2c_1^{-1}c_2^{-1}$ and $(c_1c_2)^s$ are conjugate then
$s=2$. Even if $s=2$ then two cyclically reduced words
$c_1c_2c_1^{-1}c_2^{-1}$ and $c_1c_2c_1c_2$  are cyclically
different since $c_1$ and $c_2$ are not of multiplicity $2$.

\item \textbf{$b_i+p_i=3$, $1\leq b_i\leq 2$, $S_i$ orientable and
$2\in A_i$ for all $i$ :}

This is a case where the figure eight argument does not work. The
following lemma provides us a double cover of the manifold which
has nontrivial extended products by previous cases or the non
separating torus argument. Appendix B is devoted to the proof of
this lemma.

\begin{lemma}
Let $M$ be as above, $|\mathcal{T}|=n$, $b_i+p_i=3$, $1\leq
b_i\leq 2$ and $2\in A_i$ and $S_i$ orientable for all $i$.
\begin{romlist}
\item If $n\geq3$ then $M$ has a double cover with no non
separating torus whose torus decomposition has a Seifert piece
with $3$ boundary components. \item If $n=2$ and $r\in A_1$,
$r\neq 2$, then $M$ has a double cover with no non separating
torus whose torus decomposition has a Seifert piece with two
singular fibers of multiplicity $r\neq2$. \item If $n=2$ and
$A_1=A_2=\{2\}$ then $M$ has a double cover with a non separating
torus.
\end{romlist}
\end{lemma}

\item \textbf{$n=|\mathcal{T}|$:\begin{itemize} \item $b_i=2$,
$p_i=1$, $A_i=\{2\}$ and $S_i$ is orientable for $i\neq1,n$ \item
$b_1=b_n=1$;\item $S_1$ is non orientable; \item If $S_n$ is
orientable then $p_n=1$ and $2\in A_n$
\end{itemize}}
\noindent This case is similar to the previous case except that
some $\overline{M_i}$ may not be orientable as fibration.

\noindent\textbf{\textsf{If $\overline{M_1}$ has a singular fiber
or $S_1$ has more than one cross cap:}} then consider $M'_1$ the
double cover of $\overline{M_1}$ which is Seifert manifold and is
oriented as fibration (see Remark \ref{rem-oriented-seifert}).
$M_1'$ has at least 2 singular fibers and 2 boundary components or
its base surface has genus greater thatn zero. We can construct
$\tilde{M}$, a double cover of $M$, by gluing two copies of
$M\setminus M_1$ to $M'_1$ along the boundary components of $M_1'$
such that the result double covers $M$. The covering on the
complement of $M_1'$ is the trivial double cover. Note that the
torus decomposition of $\tilde{M}$ has a Seifert piece $M_1'$
which has either two boundary components and 2 singular fibers
thus by case (3) $\tilde{M}$ has nontrivial extended loop products
or its base surface has genus greater than $0$ thus it has a non
separating torus and therefore nontrivial extended loop products.

\noindent\textsf{\textbf{If $\overline{M_1}$ has no singular fiber
and $S_1$ has one cross cap:}} then $S_1$ is indeed M\"obius strip
and $\overline{M_1}$ has another Seifert fibration model whose
base surface is a disk and has 2 singular fibers of multiplicity
2. Now if $S_n$ is orientable this is just case (5). If $S_n$ is
not orientable then it is the case above.

\end{numlist}
\end{proof}
\section{Gluing to Hyperbolic Manifolds along Torus}

In this section we analyze gluing of Seifert/hyperbolic and
hyperbolic $3$-manifolds along torus and for that we we need the
following lemma. The proof can be found in Appendix A (see
\ref{lem-hyp-group}).
\begin{lemma}
\label{lemma-ghg}
 Suppose that $G_1 ,G_2$ and $H$ are three groups with
monomorphisms $\phi_i:H\rightarrow G_i$, $i=1,2$ and $G_1*_HG_2$
denotes the amalgamated free product of $G_1$ and $G_2$ over $H$.
Let $g_1\in G_1\setminus H$ and $g_2\in G_2\setminus H$ and $
h\neq 1$ in $H$ such that:
\begin{romlist}
 \item $g_1^{-1}Hg_1\cap H=1$
\item  $g_2h\neq hg_2 $.
\end{romlist}
Then two conjugacy classes $[hg_1g_2]$ and $[hg_2g_1]$ of
$G_1*_HG_2$ are different.
\end{lemma}

\begin{proposition}
\label{hyper-hyper} Let $M$ be a connected oriented closed
$3$-manifold which contains an incompressible torus $T$ such that:
\begin{romlist}
\item $M-T$ is disconnected with two components $M_1$ and $M_2$.
\item $\overline{M_1}$ has a hyperbolic interior with finite
volume.
 \item $\overline{M_2}$ either has a hyperbolic interior with finite volume or is
 Seifert fibred not homeomorphic to $S^1 \times S^1 \times I$.
 \end{romlist}
Then $M$ has nontrivial extended loop products.
\end{proposition}
\begin{proof}
Let $\varphi:S^1\times S^1 \rightarrow T\subset M $ be a
homeomorphism. $\varphi$ gives rise to an element
$H\in\mathbb{H}_1(M)$, whose marked points are $\varphi(s,1), s\in
S^1$ and the loop passing through $\varphi(s,1)$ is
\begin{equation*}
\begin{array}{c}
\varphi_s:\S^1 \rightarrow M\\
 \varphi_s(t)=\varphi(s,t)
\end{array}
\end{equation*}

 Take $\varphi(1,1)$ as the base point for $M$, $\overline{M_1}$ and $\overline{M_2}$. Let
$h=[\varphi_1] \in \pi_1(M)$.
\begin{center}
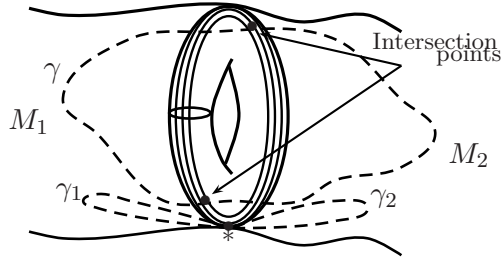
\begin{figure}[h]
\begin{pspicture}(0,-1)(15,6)
%\psgrid(0,-6)(15,6)
\psellipse[linewidth=1.2pt](6.9,2.5)(1.8,3.3)
\psellipse[linewidth=0.9pt](6.9,2.5)(1.2,3.0)
\psellipse[linewidth=0.9pt](6.9,2.5)(1.4,3.2)
\psellipse[linewidth=0.9pt](6.9,2.5)(1.6,3.3)
\pscurve[linewidth=1.2pt](7,4.2)(6.8,3.8)(6.4,2.5)(6.7,1.2)(7,0.8)
\pscurve[linewidth=1.2pt](6.9,4)(7.2,2.5)(6.8,1.1)
\pscurve[linewidth=1.2pt](1,5.7)(2,5.5)(7,5.8)(9,5.3)(11,5.5)(12,5,3)
\rput[b]{0}(0,-6.6){\pscurve[linewidth=1.0pt](1,5.7)(2,5.5)(7,5.8)(9,5.3)(11,5.5)(12,5,3)}
\psccurve[linewidth=1.0pt,linestyle=dashed](2,3)(4,4.8)(6.2,5)(10,4.8)(12,3)(13,2)(12,1)(10,0.6)(9,0)(7,-0,5)(5,0)(3,2)(2,3)
\pscurve[linewidth=1.0pt,linestyle=dashed](7,-0.8)(10,0)(11,-0.2)(10,-0.6)(7,-0.8)
\pscurve[linewidth=1.0pt,linestyle=dashed](7,-0.8)(4,0)(3,0.2)(2.6,0)(3,-0.2)(7,-0.8)
%\pscurve[linewidth=1.0pt](7,5.8)(8,2.5)(7,-0.8)
\rput[b]{0}(2.2,0){\large$\gamma_1$}
 \rput[b]{0}(11.5,-0.2){\large$\gamma_2$}
\rput[b]{0}(1.7,3.5){\Large $\gamma$}
\psellipse[linewidth=1.0pt](5.75,2.6)(0.65,0.2)
\rput[b]{0}(6.9,-1.2){$*$}
 \rput[b]{0}(7.6,5.0){\small$\bullet$}
\rput[b]{0}(6.2,-0.15){\small$\bullet$}
\rput[b]{0}(6.9,-0.9){\small$\bullet$}
\rput[b]{0}(1,2){\large$M_1$}
 \rput[b]{0}(14,1.){\large$M_2$}
\psline{->}(12,4)(7.8,4.95)
 \psline{->}(12,4)(6.4,0.2)
  \rput[b]{0}(13,4.5){\small Intersection} \rput[b]{0}(14,4) {\small points}
\end{pspicture}
\caption{$M_1$ and $M_2$ glued along $T$}\label{fig-hyp-hyp}
\end{figure}
\end{center}
We recall that $\pi_1(M)$ is the amalgamated free product of
$\pi_1(\overline{M_1})$ and $\pi_1(\overline{M_2})$ over
$\pi_1(T)$. Pick two arbitrary elements $g_1\in
\pi_1(\overline{M_1})$ and $g_2\in\pi_1(\overline{M_2})$. Let
$g_1g_2 \in\pi_1(M)$ be their product in the amalgamated free
product $\pi_1(\overline{M_1})*_{\pi_1(T)}\pi_1(\overline{M_2})$.
Let $\gamma_1$ and $\gamma_2$ be two loops in $\overline{M_1}$ and
$\overline{M_2}$ representing $g_1$ and $g_2$. The composition
$\gamma_1\gamma_2$ of these two loops is a representative for
$g_1g_2$. This loop has free homotopy type $[g_1g_2]$ (see Figure
\ref{fig-hyp-hyp}).

 Replace $\gamma_1\gamma_2$ by a loop $\gamma$ with the same free homotopy
type which intersects transversally the torus $T$ at two points
(see Figure \ref{fig-hyp-hyp}). $\gamma$ can be regarded as an
element of $\H_0(M)$. We denote this homology class by $\Gamma$.

We claim that there is a choice of $g_i \in \pi_1(\overline{M_i})
\setminus \pi_1(T)$, $i=1,2$, and $\varphi$ such that:
$$p_{A_M}(\Delta \Gamma )\bullet p_{A_M}(\Delta H)
\neq 0.$$ In the expansion of $p_{A_M}(\Delta \Gamma)\bullet
p_{A_M}(\Delta H)$, we get $4$ terms for each intersection point,
just like in the proof of Proposition \ref{selfglu-torus},
therefore $8$ terms overall. Each term is an element of $\H_0 (M)$
and represented by a loop. They have the free homotopy types:
$$
[hg_1g_2],[h],[g_1g_2],1, [hg_2g_1],[h],[g_1g_2],1.
$$

Calculating mod 2, we get $[hg_1g_2]+[hg_2g_1]$.
 By Lemma \ref{lemma-ghg}
$[hg_1g_2]$ and $[hg_2g_1]$ are different conjugacy classes if
$\pi_1(\overline{M_1})$, $\pi_1(\overline{M_2})$ and $\pi_1(T)$
satisfy the conditions of Lemma \ref{lemma-ghg}.

By assumption the interior of $\overline{M_1}$ is hyperbolic.
$\pi_1(\overline{M_1})$, can be regarded as a subgroup of
$PSL(2,\mathbb{C})$ where the image of $\pi_1(T)$ in
$\pi_1(\overline{M_1})$ consists of parabolic elements which have
a common fixed point.

 For $g\in \pi_1(\overline{M_1})\setminus \pi_1(T)$ (we have denoted the image of
$\pi_1(T)$ in $\pi_1(\overline{M_1})$ simply by $\pi_1(T)$),
$g^{-1} \pi_1(T) g \cap \pi_1(T)=\{1\}$, since if the fixed point
of the elements of $\pi_1(T)$ is $m$ then the elements of $g^{-1}
\pi_1(T) g$ fix $g^{-1}(m)$. So if $g^{-1} \pi_1(T) g \cap
\pi_1(T)$ is nontrivial then $m$ must be a fixed point of $g$ but
this implies that $g\in \pi_1(T)$ which is a contradiction. Thus
$G_1=\pi_1(M_1)$ satisfies condition $(i)$ of Lemma
\ref{lemma-ghg}.

Also, the argument above proves that if $\overline{M_2}$ has
hyperbolic interior then for any $g_2\in
\pi_1(\overline{M_2})\setminus \pi_1(T)$ and $h\in
\pi_1(T)\setminus \{1\}$ the condition $(ii)$ is satisfied. So the
assertion is proved if $\overline{M_1}$ and $\overline{M_2}$ both
have hyperbolic interior with finite volume.

Now consider the case when $\overline{M_2}$ is Seifert fibred. If
$h$ is not an power of the element by represented by the normal
fiber of $\overline{M_2}$, then it is not in the centralizer of
$\pi_1(\overline{M_2})$ and the condition $(ii)$ holds for an
appropriate choice of $g_2$. Otherwise replace $\varphi$ with
$$
\varphi'(s,t)=\varphi(t,s).
$$

Namely, the set of marked points is:
 $$\varphi'(s,1)=\varphi(1,s), 0\leq s\leq 1$$ and the loop passing through
 $\varphi'(s,1)$ is $$ \varphi'_s(t)=\varphi'(s,t)=\varphi(t,s), 0\leq t \leq 1. $$

Then $h$, the element of $\pi_1(\overline{M_2})$ represented by
loop passing through $\varphi'(1,1)$, is not in the centralizer of
$\pi_1(\overline{M_2})$ and one can choose  $g_2\in
\pi_1(\overline{M_2})\setminus\pi_1(T)$ such that $g_2h\neq hg_2$.
This finishes the proof.
\end{proof}
\section{Proof of Main Theorem: Part I}

In this section we prove the first part of the main theorem,
namely that if a closed oriented manifold $M$ which is not
algebraically hyperbolic, then $M$ or a double cover of $M$ has
nontrivial extended loop products.

For this we use the \textit{prime decomposition} and the
\textit{theory of characteristic surfaces} for 3-manifolds. We
recall some of these theories.

A closed oriented $3$-manifold $M$ is said to be \emph{prime} if
it cannot be obtained as the connected sum of two $3$-manifolds
that none of them is homeomorphic to $3$-sphere. By the
Alexander's theorem $S^3$ is prime, as every embedded $2$-sphere
in $S^3$ bounds a 3-ball.

\textbf{Prime Decomposition} (see \cite{Mil}): \label{theorem3}
Let $M$ be a compact, connected and oriented $3$-manifold. Then
there is a decomposition $$M=P_1\#P_2\#...\#P_n$$ with each $P_i$
prime and this decomposition is unique up to permutation of
$P_i$'s and insertion or deletion of $S^3$'s.

A closed $3$-manifold is said to be \emph{irreducible} if any
embedded $2$-sphere $S^2$ bounds a $3$-ball. Again by the
Alexander's theorem, $S^3$ is irreducible. An oriented irreducible
$3$-manifold is prime. The converse is not true. The only prime
oriented $3$-manifold which is not irreducible, is $S^1 \times
S^2$.

One can deduce a necessary condition for irreducibility from the
\emph{sphere theorem}. Here, we state a simplified version of the
theorem.

\textbf{Sphere Theorem} (see [6]): Let $M$ be an oriented
connected $3$-manifold. If $\pi_2(M)\neq0$ then there is an
embedded $S^2$ in $M$ representing a nontrivial element in
$\pi_2(M)$.

The topology of an irreducible $3$-manifold $M$ is coarsely
determined by the cardinality of the fundamental group. By the
\emph{sphere theorem} (see \cite{He}) $\pi_2(M)=0$. Let
$\tilde{M}$ be the universal cover of $M$. If $\pi_1(M)$ is
finite, then $\tilde{M}$ is closed and simply connected. Therefore
$\tilde{M}$ is a homotopy sphere. If $\pi_1(M)$ is infinite then
$\tilde{M}$ is not closed and we have $\pi_i(\tilde{M})=0 $ for
all $i$, by Hurewitcz theorem. Therefore $\tilde{M}$ is
contractible or in other words $M$ is a $K(\pi,1)$(aspherical).

Now by the prime decomposition theorem, every $3$-manifold $M$ can
be written as
\begin{equation}\label{prime-decomposition}
M= (K_1\#K_2\#...\#K_p)\#(L_1\#L_2...\#L_q) \#(\#_1^r (S^1\times
S^2))
\end{equation}
where $K_i$'s are closed aspherical irreducible $3$-manifolds with
infinite $\pi_1$ and $L_i$ are closed $3$-manifolds that are
finitely covered by a homotopy $3$-sphere.

 Let $S$ be an embedded
surface in a $3$-manifold $M$. We shall say that $S$ is
\emph{properly embedded} if $S \cap \partial M=\partial S$. Then
$S$ is \emph{two-sided } in $M$ if the normal bundle of $S$ in $M$
is trivial.

\noindent A $2$ sided surface $S$ in $M$ is said to be
\emph{compressible} if either:
\begin{romlist}
\item $S$ is a $2$-sphere and it bounds a 3-cell in  $M$\\
or \item there is a disk $D \in M$  such that $D \cap S=\partial
D$ and $\partial D$ is non contractible in $S$.
\end{romlist}
Otherwise $S$ is said to be \emph{incompressible}. An embedded
disk is considered to be incompressible. For an incompressible
surface $S$ in $M$, the induced map by inclusion
$\pi_1(S)\rightarrow \pi_1(M)$ is injective. Indeed, a $2$-sided
surface $S$ is incompressible if and only if the map
$\pi_1(S)\rightarrow \pi_1(M)$ induced by inclusion is
injective\footnote{$\pi_1$-injective}.

 An oriented $3$-manifold $M$, possibly with boundary is called
\emph{atoroidal} if any incompressible torus in $M$ is
\emph{$\partial$-parallel}, meaning that it can be isotoped to the
boundary of $M$.

 \textbf{Torus decomposition} (see \cite{JS} or \cite{Johannson}):
For a compact irreducible oriented $3$-manifold $M$ there exists a
collection $\mathcal{T} \subset M$ of disjoint incompressible tori
such that the closure of each component of $M \backslash
\mathcal{T}$ is either atoroidal or a Seifert manifold with
$\mathcal{T}$ tangent to the fibration. A minimal such a
collection $\mathcal{T}$ is unique up to isotopy.

As we are interested in rank 2 abelian subgroups of fundamental
groups, torus theorems by Scott (\cite{PScott}) and Seifert fibred
space theorem by Casson-Jungeris \cite{CJ} and Gabai \cite{Gab}
will be useful to us.

\textbf{Torus Theorem} (see \cite{PScott}): Let $M$ be a compact,
connected oriented, irreducible $3$-manifold. Suppose that
$\pi_1(M)$ has a free abelian subgroup of rank $2$. Then either
\begin{romlist}
\item $M$ contains a two sided torus such that the inclusion
homomorphism $\pi_1(T)\rightarrow \pi_1(M)$ is injective, or \item
$\pi_1(M)$ contains an infinite cyclic normal subgroup.
\end{romlist}
\textbf{Seifert Fibred Space Theorem} (see \cite{CJ} or
\cite{Gab}) Let $M$ be a compact connected, oriented, irreducible
$3$-manifold with infinite $\pi_1$, then $M$ is a Seifert fibred
space if and only if $\pi_1(M)$ has a infinite cyclic normal
subgroup.

We recall the definition of algebraically hyperbolic manifolds.

\noindent\textbf{Definition} An oriented 3-manifold is said to be
\textit{algebraically hyperbolic} if it is aspherical ($K(\pi,1)$)
and $\pi_1(M)$ has no rank 2 abelian subgroup.\\

\noindent\emph{Proof of the first part of the main theorem:}
 If $M$ is not algebraically hyperbolic then at least one of
the following holds:

\subsection{If \textit{M} is not aspherical}\label{sec-not-aspherical}
By prime decomposition theorem, every $3$-manifold $M$ can be
written as
\begin{equation}\label{prime-decomposition}
M= (K_1\#K_2\#...\#K_p)\#(L_1\#L_2...\#L_q) \#(\#_1^r (S^1\times
S^2))
\end{equation}
where $K_i$'s are closed aspherical irreducible $3$-manifolds with
infinite $\pi_1$ and $L_i$ are closed $3$-manifolds that are
finitely covered by a homotopy $3$-sphere. If $M$ is not
aspherical then one of the following holds:
\begin{romlist}
\item If $M$ is a connected sum of manifolds with nontrivial
fundamental group then Proposition \ref{connected-sum} says that
$M$ has nontrivial extended loop products.

 \item If $M$ has a non separating 2-sphere or in other words there is a $S^2\times S^1$ factor in its
prime decomposition, then by Corollary \ref{nonseparating2} $M$
has nontrivial extended loop products.

\item $M$ has finite fundamental group; Then by Proposition
\ref{finitegroup} $M$ has nontrivial extended loop products.

\end{romlist}
Thus in either cases  $M$ has nontrivial extended loop product.

\subsection{If $M$ aspherical and $\Z\oplus\Z\subset \pi_1(M)$}

We turn to the case when $M$ is aspherical $3$-manifolds with rank
$2$ abelian groups. If follows from the argument in
\ref{sec-not-aspherical} that $M$ is either irreducible or or is
connected sum of an irreducible  3-manifold and homotopy
3-sphere\footnote{ Note that such a connected sum is again
aspherical.}.

\subsubsection{M is irreducible}
\label{case-irreducible}

By Torus theorem and Seifert theorem one of the following holds:

\subsection*{\emph{a) M is a closed Seifert manifold}}
 If $M$ is Seifert manifold then by Proposition \ref{seifert-closed} it has nontrivial extended loop products.

\subsection*{\textit{b) The collection of tori in the torus
decomposition of $M$ is nonempty.}}

This case itself is divided into $3$ cases:

\subsection*{1.  \textit{M}\textsf{ contains a non separating torus:}}

If there is a non separating torus in the collection $T$ then by
Proposition \ref{selfglu-torus} $M$ has nontrivial extended loop
products.

 \subsection*{\textsf{2. Torus decomposition of} \textit{M} \textsf{has
only Seifert pieces and has no non separating torus:}}

Let $\mathcal{T}$ be the collection of tori provided by the torus
decomposition of $M$ and
$$
M\setminus\mathcal{T}=M_1\cup M_2...\cup M_n.
$$
$\overline{M_i}$ is a Seifert manifold. Let $p_i$ and $b_i$  be
 the number of singular fibers and boundary components of
$\overline{M_i}$. If $S_i$, the base surface of $\overline{M_i}$,
is orientable then $g_i$ denotes its genus otherwise $g_i$ is
number of cross caps in $S_i$.

If for some $i$, $S_i$ is orientable and $g_i\geq 1$, then $M$ has
a non separating torus and by case 1 we know that $M$ has
nontrivial extended loop product. So we assume that $g_i=0$ if
$S_i$ is oriented.

 If $S_i$ is
orientable then,
 $$ p_i+b_i\geq 3.
$$
 To see that suppose that $p_i+b_i\leq 2$. Then one of the
 following holds:

 \begin{alphlist}
\item $b_i=1$ and $p_i=1$; \item $b_i=2$ and $p_i=0$
\end{alphlist}

In the first case if follows that $\pi_1(\overline{M_i})\simeq \Z$
so the boundary components of $\overline{M_i}$ can not be
incompressible. In the latter case it follows that
$\overline{M_i}\simeq S^1\times S^1 \times [0,1]$  this
contradicts the minimality of the collection $\mathcal{T}$ as one
can extend the Seifert fibration of one of the neighboring
$\overline{M_j}$ to $\overline{M_i}$.

Therefore $M$ satisfies all the conditions of Proposition
\ref{seifert-glued}, hence $M$ or a double cover of $M$ has
nontrivial extended loop product.

\subsection*{3. \textsf{If} \textit{M} \textsf{ has an atoroidal component and has no non separating torus:}}
 Suppose that $\mathcal{T}$ is the collection of tori
provided by the torus decomposition of $M$ and every torus in
$\mathcal{T}$ is separating.

By Thurston theorem \footnote{Geometrization of Haken
manifolds.}\cite{T},

\noindent \textbf{Theorem} \emph{A compact, oriented, irreducible
atoroidal 3-manifold whose boundary consists of a finite number of
tori admits a complete hyperbolic metric of finite volume;}\\
atoroidal components of $M\setminus \mathcal{T}$ are hyperbolic.

Let $M_1$ be a hyperbolic component of $M\setminus \mathcal{T}$
and $M_2$ is another component
 that is attached to $M_1$ along a torus $T$. By proposition
 \ref{hyper-hyper}, $M'=\overline{M_1}\cup T_1\cup \overline{M_2}$ has nontrivial extended loop products.

 Consider the same homology classes $\Gamma $ and $H$ in $\H _* (M')$ constructed in
 the proof of Proposition \ref{hyper-hyper}. They can be regarded as
 elements of $\H _* (M)$. We claim that the loop product

\begin{equation*}p_{A_M}(\Delta(\Gamma ))\bullet p_{A_M}(\Delta H)
\neq 0,\end{equation*} so $M$ has nontrivial extended loop
products.

The calculation of the loop product is the same as of the proof
proposition  \ref{hyper-hyper}. Calculating mod 2 only two terms
remain, each one an element of $\H_0(M)$ which is represented by a
conjugacy class in $\pi_1(M)$. These conjugacy classes are
\begin{equation*}
 [hg_1g_2] \text{ and } [hg_2g_1].
\end{equation*}

We prove that these two conjugacy classes are different. By the
proof of proposition  \ref{hyper-hyper} we know that $hg_1g_2$ and
$hg_2g_1$ are not conjugate in $\pi_1(M')$. So by the injectivity
lemma (Lemma \ref{injection-lem}), if we show that $hg_1g_2$ is
not conjugate to an element of the fundamental group of one of the
boundary components of $\overline{M'}$ it follows that $hg_1g_2$
and $hg_2g_1$ are not conjugate in $\pi_1(M)$.

Let $T_1\cup T_2...\cup T_k=\partial \overline{M'}$ and suppose
that $hg_1g_2$ is conjugate to an element of $\pi_1(T_i)$, $i>1$
in $\pi_1(M')$. Since one of the
\begin{equation*}
T_i\subset \overline{M_1} \text{ or } T_i\subset \overline{M_2}
\end{equation*}
holds let us assume that $T_i\subset \overline{M_1}$. Then
\begin{equation*}
\pi_1(T_i)\subset \pi_1(M_1)= \pi_1(\overline{M_1})
\end{equation*}
as $T_i$ is incompressible.  If $hg_1g_2$ is conjugate to an
element to $\pi_1(T_i)$ then it is conjugate to an element in
$\pi_1(M_1)$. But this is a contradiction since
$hg_1g_2=(hg_1)g_2$ as an element of
\begin{equation*}
\pi_1(M')=\pi_1(M_1)*_{\pi_1(T)}\pi_1(M_2)
\end{equation*}
is a cyclically reduced word of length 2 and it is cannot be
conjugate to an element of the factor $\pi_1(M_1)$. Now by the
injectivity lemma $hg_1g_2$ and $hg_2g_1$ are not conjugate in
$\pi_1(M)$.

\subsubsection{$M$ is a connected sum of an irreducible
aspherical manifold and a homotopy 3-sphere}

Suppose that $M=N\#P$ where $N$ is irreducible and aspherical and
$P\neq S^3$ is a homotopy 3-sphere. $M$ is obtained from $N$ by
replacing a standard ball $D\subset N$ with a fake ball $D'$
namely $M=(N\setminus D)\cup_{S^2}D'$. Let
$$
f: N \rightarrow M
$$
be a homotopy equivalence which sends $D$ to $D'$ and is identity
on the complement of $D$ in $N$. Since $\Z\oplus
\Z\subset\pi_1(M)$ hence $\Z\oplus \Z\subset \pi_1(N)$.

By Seifert and Torus theorem either $N$ is Seifert or the
collection of tori in its torus decomposition is non empty.

If $N$ has a nontrivial torus decomposition it follows from the
previous case, \ref{case-irreducible}, that $N$ has nontrivial
extended loop products. The homology classes in $N$ constructed in
\ref{case-irreducible} can be regarded as homology classes in $M$
since they can have chain representative which is in the
complement of $D$ and the loop products calculated there are also
nontrivial in $M$ as $\pi_1(M)\simeq \pi_1(N)$.

If $N$ is a Seifert manifold then one of the 3-homology class
constructed in according to Section \ref{section-closed-seifert},
$\theta$, can not be considered as a homology class in $M$ since
it uses all the manifold $N$ but we claim that
$$
p_{A_M}(f_{\L}(\theta) )\bullet p_{A_M}( f_{\L}(\alpha))\neq 0.
$$
which proves that $M$ has nontrivial extended loop products.

 The proof is exactly the one in Section \ref{section-closed-seifert}
one just has to notice that we can choose a representative for
$\alpha$ such that its marked point is in the complement of $D$.

\section{Proof of Main Theorem: Part II, Algebraically Hyperbolic Manifolds}

\label{section-homoto-hyperbolic}
\medskip
In this section we prove the second part of the main theorem,
namely that algebraically hyperbolic manifolds have trivial
extended loop products.

 It follows from Lemma \ref{lem-compoenents} that

\begin{lemma}
\label{lemma3} If $M$ is a $K(\pi,1)$ (aspherical) then
$\mathbb{L}M$ is homotopy equivalent to
\[
\coprod_{[\alpha]\neq [1] \in \hat{\pi}_1(M)}K(\pi_{\alpha},1) \]
where we choose only one representative for each conjugacy class
and $\pi_{\alpha}$ is the centralizer of $\alpha$ in $\pi_1(M)$.
\end{lemma}
Therefore to prove the second part of the theorem we need a good
understanding of the centralizers. For that we need following
theorem.

 \textbf{Compact Realization Theorem} (see
\cite{He}): Let $M$ be a connected, oriented $3$-manifold, and let
$K$ be a subgroup of $\pi_1(M)$ which is finitely generated. Then
there is a compact, oriented irreducible $3$-manifold $M_0$ with
$\pi_1(M_0)\simeq K$.
\begin{lemma}
\label{lemma4}
 If a closed oriented $3$-manifold $M$ is algebraically hyperbolic then the fundamental group of $M$ does not have any
nontrivial finite subgroups and the centralizer $\pi_{\alpha}$ of
an element $\alpha \in \pi_1(M)$ is isomorphic to an additive
subgroup of $\mathbb{Q}$.
\end{lemma}
\begin{proof}
Let $\mathbb{Z}_p$ be a nontrivial finite subgroup of  $\pi_1(M)$.
Then $K(\mathbb{Z}_p,1)$ is a covering of $M$, so it has to be a
$3$-manifold. But one knows that $K(\mathbb{Z}_p,1)$ has homology
in infinitely many dimensions. So $\pi_1(M)$ has no finite
subgroup.

For the second part we first prove that every finitely generated
subgroup of $\pi_{\alpha}$ is isomorphic to $\mathbb{Z}$. Suppose
$\beta_1, \beta_2,....\beta_n$ are elements of $\pi_{\alpha}$.
Consider the subgroup $\langle \beta_1, \beta_2,....\beta_n
,\alpha \rangle$ which is finitely generated and it is not a free
product since it has a nontrivial center. Therefore by compact
realization theorem it is isomorphic to the fundamental group of a
compact, oriented irreducible $3$-manifold. The cyclic subgroup
$\langle \alpha \rangle $ is in its center and therefore it is
normal. Thus by the Seifert fibred space theorem, the compact
realization should be a Seifert manifold. Since $\pi_1(M)$ has no
finite subgroup or rank 2 abelian subgroup this Seifert manifold
has to be a solid torus \footnote{$M$ cannot be $S^2\times S^1$
since it is irreducible.}. As the fundamental group of a solid
torus is $\mathbb{Z}$, we have $\langle \beta_1,
\beta_2,....\beta_n ,\alpha \rangle \simeq \mathbb{Z}$. To prove
that $\pi_{\alpha}$ is isomorphic to a subgroup of $\mathbb{Q}$,
note that since $\pi_1(M)$ is countable and every finitely
generated subgroup is cyclic, we can write $\pi_{\alpha}=G_1 \cup
G_2 ....$ where each $G_i$ is an infinite cyclic subgroup and
$G_i\subset G_{i+1}$. Let $x_i$ be a generator for $G_i$; then we
have $x_i=x_{i+1}^{n_{i}}$ for some $n_i$. We construct a map $
\varphi: \pi_{\alpha} \rightarrow \mathbb{Q}$ by letting
$\varphi(z_1)=1$ and $\varphi(z_{i+1})=\frac{1}{n_1n_2...n_i}$.
This is a well defined map because of the way we chose $n_i$'s.
$\varphi$ is injective. To see that note that if $x\in \pi_1(M)$
then $x\in G_i$ for some $i$. So $x=kz_i$ for some $i$ and $k$.
$\varphi(x)=\frac{k}{n_1n_2...n_i}$. Therefore $\varphi(x)=0$
implies $k=0$ or $x=0$.
\end{proof}

\begin{lemma}
\label{lemma5} If $G$ is an additive subgroup of $\mathbb{Q}$ then
$K(G,1)$ has homological dimension one.
\end{lemma}
\begin{proof}One can write $G=\underset{n \rightarrow \infty}{\lim}G_n$, where $G_i$'s are cyclic
subgroups of $G$. So $K(G,1)=\underset{n\rightarrow
\infty}{\lim}K(G_n,1)$. Since each $K(G_n,1)\simeq K(\Z,1)$ has
homological dimension $1$  therefore  the direct limit $K(G,1)$
has homological dimension one as taking homology  commutes with
taking direct limit.
\end{proof}
\begin{proposition}
\label{homotopy-hyper-reduced} If $M$ is algebraically hyperbolic
then any finite cover of $M$ has trivial extended loop products.
\end{proposition}
\begin{proof}
Let $\widetilde{M}$ be a finite cover of $M$. Note that
$\widetilde{M}$ is also algebraically hyperbolic.
 Since $\widetilde{M}$ is a $K(\pi',1)$ then its free loop space is homotopy equivalent to
$$ \coprod_{[\alpha] \in \hat{\pi'}}K(\pi'_{\alpha},1)$$
where we choose only one representative for each conjugacy and
$\pi'_{\alpha}$ is the centralizer of $\alpha$ in
$\pi_1(\widetilde{M})$. So
\[
A_{\tilde{M}}\simeq H_*(\coprod_{[\alpha] \in \hat{\pi'},
\alpha\neq 1 }K(\pi'_{\alpha},1))\simeq \bigoplus_{[\alpha] \in
\hat{\pi'}, \alpha\neq 1 } H_*(K(\pi'_{\alpha},1))
\]
 Now by Lemma \ref{lemma4} and Lemma \ref{lemma5}, we see that each
$K(\pi'_{\alpha},1)$, $ \alpha \neq 1$, has homological dimension
$1$. This proves that the loop product $\bullet$ on
$A_{\tilde{M}}$ is zero because of the absence of homology classes
of sufficiently high degree.
\end{proof}

\appendix

\section{Some Facts in Free Products with Amalgamation}
\subsection{Some Basic Facts about Free Products with
Amalgamation:} We recall some definitions and theorems from group
theory related to the free product with amalgamation. We follow
the terminology and notations of \cite{MKS} and \cite{LS}.

Let $G_1$, $G_2$ and $H$ be three groups and $\phi_i:H\rightarrow
G_i$, $i=1,2$, be two injective homomorphisms. We can regard $H$
as a subgroup of each $G_i$, $i=1,2$ and we identify each element
of $H$ with its image under $\phi_i$ while we are in $G_i$. So if
$h \in H$, then $h$ can be considered as an element of both $G_1$
and $G_2$. Following \cite{LS}, we can form the amalgamated free
product $G_1*_HG_2$ of $G_1$ and $G_2$ over $H$.

\noindent A sequence of elements $g_1,g_2,...,g_n \in G_1*_HG_2$
is called \textit{reduced} if:
\begin{numlist}
\item Each $g_i$ is in one of the factors $G_1$ or $G_2$;
 \item $g_i,g_{i+1}$ come from different factors; So
 if $n>1$ then $g_i \notin H$ for all $i$.

\item If $n=1$ then $g_1\neq 1$.
\end{numlist}

 $g_1g_2\cdots g_n\in G_1*_HG_2$ is called a \emph{reduced word
 } if $g_1,g_2,...,g_n$ is a reduced sequence.

There is no canonical way of presenting the elements of
$G_1*_HG_2$ uniquely. Once a choice of left coset representatives
for $G_1/H$ and $G_2/H$ is made then every element can represented
uniquely. Let $c_1,....c_n$ be a right coset representative for
$G_1/H$ and $G_2/H$. Here is an informal description of the
representation:

 Let $g$ be an element in $G_1*_HG_2$ and $g=g_1g_2..g_r$ be a
 reduced presentation of $g$. Start with $g_r$, there is an $h_r
 \in H$ and $c_r \in G_1/H$ or $G_2/H$ (depending in which factor is $g_r$) such that

$$
g_r=h_rc_r.
$$
Now by replacing $g_r$ by $h_rc_r$ we have
$g=g_1g_2...g_{r-1}h_rc_r$.  Again there is a $h_{r-1}\in H$ and
$c_{r-1}$ such that
$$
g_{r-1}h_r=h_{r-1}c_{r-1}.
$$
By replacing $g_{r-1}h_r$ by $h_{r-1}c_{r-1}$ we have
$g=g_1g_2...g_{r-2}h_{r-1}c_{r-1}c_r$. Continuing this procedure
we  end up with $g=h_1c_1c_2...c_r$. This is the desired
presentation. While stating it more precisely, the following
theorem also asserts its uniqueness (\emph{i.e.} its independence
from the reduced word we used in the construction).
\begin{theorem}
(See \cite{MKS} pp.201) Let $G=G_1*_HG_2$ where $H$ is a mutual
subgroup of $G_1$ and $G_2$. Suppose a specific right coset
representative systems for $G_1/H$ and $G_2/H$ have been selected.
Then with each element of $g$ of $G$ we can associate a unique
sequence\footnote{This presentation is called \emph{reduced form}
(see \cite{MKS}).}\end{theorem}
\begin{equation*}
(h,c_1,c_2, ....c_r)
\end{equation*}
such that:
\begin{romlist}
\item  $h$ is an element, possibly 1, of $H$; \item $c_i$ is a
coset representative of $G_1/H$ or $G_2/H$; \item $c_i\neq 1 $
\item $c_i$ and $c_{i+1}$ are not both in $G_1$ or $G_2$. \item If
$h'$ and $c_i'$s are the images of $h$ or $c_i$'s under the
homomorphisms of $G_1$ or $G_2$ into $G$, then
$g=h'c_1'c_2'...c_r'$ in $G$.
\end{romlist}

 A sequence of elements $g_1,g_2,...,g_n \in G_1*_HG_2$ is called
\textit{cyclically reduced} if all cyclic permutations of the
sequence $ g_1,g_2,...,g_n$ are also reduced.

 $g_1g_2\cdots g_n\in G_1*_HG_2$ is called a \emph{cyclically reduced word
 } if $g_1,g_2,...,g_n$ is a cyclically reduced sequence.

Also every element of $G_1*_HG_2$ is conjugate to $g_1g_2\cdots
g_n$ where the sequence $g_1,g_2,...,g_n$ is cyclically reduced.
Indeed if $g=g'_1g'_2\cdots g'_k$  where $g'_1,g'_2,...,g'_k$ is a
reduced sequence. If $g'_1 $ and $g'_k$ are from
 different factors, then clearly $g'_1g'_2\cdots g'_k$ is cyclically reduced.
 Otherwise consider $g'_kg{g'_k}^{-1}=g'_kg'_1g'_2\cdots g'_{k-1}$ and replace
 $g'_kg_1$ by $g'=g'_kg'_1$ which is from the same factor. Now
 $g',g'_2,...,g'_{k-1}$ is cyclically reduced.

\begin{theorem}(see \cite{LS} p.187) If $c_1,c_2,...c_n$, $n\geq 1$, is a
reduced sequence in $G=G_1*_HG_2$, then the product $c_1c_2\cdots
c_n\neq 1$ in $G$.
\end{theorem}

The following theorem describes the conjugacy classes of
$G_1*_HG_2$.

\begin{theorem} (see \cite{MKS} pp.212)

 Let $G=G_1*_HG_2$. Then every element of $G$ is conjugate
to a cyclically reduced element of $G$. Moreover, suppose that $g$
is a cyclically reduced element of $G$. Then:
\begin{romlist}
\item If $g$ is conjugate to a cyclically reduced word
$p_1p_2..p_r$ where $r\geq 2$; then $g$ can be obtained by
cyclically permuting $p_1,p_2,...,p_r$ and then conjugating by an
element of $H$. \item If $g$ is conjugate to an element $h$ in
$H$, then $g$ is in some factor and there is a sequence of
$h,h_1,...,h_t,g $ where $h_i$ is in $H$ and consecutive terms of
the sequence are conjugate in some factor. \item If $g$ is
conjugate to an element $g'$ in some factor but not in a conjugate
of $H$ then $g$ and $g'$ are in the same factor and are conjugate
in that factor.
\end{romlist}
\end{theorem}
\begin{lemma}
\label{lem-hyp-group} Suppose that $G_1 ,G_2$ and $H$ are three
groups with monomorphisms $\phi_i:H\rightarrow G_i$, $i=1,2$. Let
$G_1*_HG_2$ denote the amalgamated free product of $G_1$ and $G_2$
over $H$. Let $g_1\in G_1\setminus H$ and $g_2\in G_2\setminus H$
and $ h\neq 1$ in $H$ such that:
\begin{romlist}
 \item $g_1^{-1}Hg_1\cap H=1$
\item  $g_2h\neq hg_2 $.
\end{romlist}
Then two conjugacy classes $[hg_1g_2]$ and $[hg_2g_1]$ of
$G_1*_HG_2$ are different.
\end{lemma}
\begin{proof}
Suppose that $[hg_1g_2]=[hg_2g_1]$. It is clear that
\begin{align}
[hg_1g_2]=[g_2hg_1] \\
[hg_2g_1]=[g_2g_1h]
\end{align}

 So $g_2(hg_1)$ is conjugate to $g_2(g_1h)$. Since $g_2(g_1h)$ is cyclically reduced  (of
 length $2$, $p_1=g_2$ \& $p_2=g_1h$) by theorem A.3 part $(ii)$,
$g_2(g_1h)$ is conjugate to $g_2(g_1h)$ or $(g_1h)g_2$ by an
element in $H$.

In the latter case we have
$$
g_2(hg_1)=h_1(g_1h)g_2h_1^{-1}
$$
which is equivalent to
$$
g_2^{-1}(h_1g_1h)^{-1}g_2(hg_1h_1)=1.
$$

But this contradicts Theorem A.2. Since the sequence
$$g_2^{-1},(h_1g_1h)^{-1},g_2,(hg_1h_1)$$ is a reduced sequence.

Therefore there should be an element $h_1\in H$ such that
\begin{equation*}
g_2hg_1=h_1g_2g_1hh_1^{-1}
\end{equation*}
 or equivalently

\begin{equation*}
 g_2hg_1h_1=h_1g_2g_1h
\end{equation*}
or
\begin{equation}
\label{eqconjugacy}
 (g_2h)(g_1h_1)=(h_1g_2)(g_1h)
\end{equation}

Now using the uniqueness of representation in terms of left cosets
of $H$ Theorem 1 (also see the construction) the equation
(\ref{eqconjugacy}) implies that $g_1h_1$ and $g_1h$ are in the
same left coset in $G_1/H$ or in other words, there exists $h_2\in
H$ such that
\begin{equation*}
g_1h_1=h_2g_1h
\end{equation*}
or equivalently
\begin{equation*}
h_1h^{-1}=g_1^{-1}h_2g_1
\end{equation*}

So $h_1h^{-1}=g_1^{-1}h_2g_1 \in g_1^{-1}Hg_1 \cap H=1$ hence
$h_2=h_1h^{-1}=1$ which implies $h_1=h$. By substituting  $h_1=h$
in (\ref{eqconjugacy}) and simplifying, it follows
\begin{equation*}
g_2h=hg_2
\end{equation*}
 which is a contradiction. Therefore $[hg_1g_2] \neq [hg_2g_1]$

\end{proof}
\subsection{Tree of Groups}

A \textit{tree of groups} $(G,T)$ consists of a finite tree $T$
and $G$ a collection of groups
$$
G=\{G_e\}_{e\in \text{edge} T}\coprod \{G_v\}_{v\in \text{vert}
T},
$$
 a group $G_v$ for every vertex $v\in $ vert $T$, a group $G_e$
for every edge $e$ in $\edg T$ and a monomorphism
$$f_e^v:G_e\rightarrow G_v$$
 if $v$ is a vertex of the edge $e$ and

We call $G_v$ a \textit{vertex group} if $v$ is a vertex of $T$
and $G_e$ a \textit{edge group} if $e$ is an edge of $T$. If $v$
is a vertex of the edge $e$  then $G_e$ can be considered as a
subgroup of $G_v$.

Suppose that $(G,T)$ is a tree of groups, $v_0$ a vertex of $T$.
Let $n=|\edg T|$. Consider a sequence of trees $T_i$, $1\leq i\leq
n$ such that
\begin{romlist}
\item $|\edg T_i |=i$, $0\leq i\leq n$
 \item $\vert T_0=\{v_0\}$,
 \item $T_i\subset T_{i+1}$ and
 \item $T_n=T$.

 \end{romlist}
Let $e_1,e_2,...e_n$ be the sequences of the edges and
$v_1,v_2,...v_n$ be the sequence of vertices such that
\begin{itemize}

 \item $\edg T_1=\{e_1\}$
 \item $\vert T_1=\{v_0,v_1\}$

\item $\edg T_{i+1}=\edg T_i\cup \{e_{i+1}\}, \quad 1\leq i \leq
n-1. $ \item $\vert T_{i+1}=\vert T_i\cup \{v_{i+1}\} \quad 1\leq
i \leq n-1$.
\end{itemize}
\noindent Let
$$
G_{T_0}=G_{v_0},
$$
\begin{center}
and\\

$ G_{T_i}=G_{T_{i-1}}*_{G_{e_n}}G_{v_n} $

$ 1 \leq i \leq n$
\end{center}

$G_T=G_{T_n}$ is called  the amalgamation of $G_v$'s along $G_e$'s
and is independent of the choice of the sequence $T_i$ and depends
only on $(G,T)$ and there is an inclusion

$$
\varphi_k:G_k\rightarrow G_T
$$
 where $k$ is a vertex or an edge.\\
The following lemma is a generalization of theorem A.3 part (iii).
\begin{lemma}
\label{lemma-tree}
 Let $(G,T)$ be a tree of groups and
 $v_0$ a vertex of $T$. Consider the construction above.
Suppose that $[a]$ and $[b]$ are two distinct conjugacy classes of
$G_{v_0}$ and $a$ is not conjugate in $G_{v_0}$ to an element of
any edge group $G_e\subseteq  G_{v_0}$ where $v_0$ is a vertex of
$e$. Then $[a]$ and $[b]$ as the conjugacy classes of $G_T$ are
distinct.
\end{lemma}
\begin{proof}
$T_i$'s, $v_i$'s and $e_i$'s are as above.

 We prove the lemma by
induction. We prove that for every $i$, $1 \leq i \leq n$

\begin{numlist}
\item $[a]$ and $[b]$ are distinct as conjugacy classes of
$G_{T_i}$.

\item $a$ is not conjugate in $G_{T_i}$ to an element in the
vertex group $G_w$, where $w \neq v_0$ is a vertex of $T_i$.
\end{numlist}

We verify the statement for $i=1$. By the assumption and part
$(iii)$ of theorem A.3 it follows that the first statement of
induction is true for $i=1$. For the second part, suppose that $a$
is conjugate in $G_{T_1}$ to an element of $G_{v_1}$, then again
by the third part of theorem A.3, $a$ has to be conjugate to an
element of $G_{e_1}$ and then by part $(ii)$ it has to be
conjugate in $G_{v_0}$ to an element of the vertex group $G_{e_1}$
which contradicts our assumption.

Suppose that the statement is true for $i$. By the second
statement of the induction for $i$, $a$ in not conjugate in
$G_{T_i}$ to an element of $G_{e_{i+1}}\subseteq G_{T_i}$.
Therefore by part $(iii)$, $a$ and $b$ represent different
conjugacy classes in $G_{T_{i+1}}$ as
$$
G_{T_{i+1}}=G_{T_i}*_{G_{e_{i+1}}}G_{v_{i+1}},
$$
proving the first statement of the induction for $i+1$. To prove
the second statement for $i+1$, suppose that by contrary $a$ is
conjugate to an element of $G_w$, where $w\neq v_0$ and $w\in
T_{i+1}$.

Since
$$
G_{T_{i+1}}=G_{T_i}*_{G_{e_{i+1}}}G_{v_{i+1}},
$$
and $a$ is not conjugate in $G_{T_i}$ to an element of
$G_{e_{i+1}}$ then again by part $(iii)$ of theorem 3 we must have
$w\in T_i$ and $a$ is conjugate in $G_{T_i}$ to an element of
$G_w$. This contradicts the assumption of the induction for $i$.
Hence the second statement of the induction holds for $i+1$.\\
Finally for $i=n$ we get the statement of the lemma.
\end{proof}

\section{Finite Covers of Certain Seifert Manifolds with Boundary}
 This appendix is devoted to explaining some facts that are used in the
proof of main theorem, section $9$. We follow the notations
introduced section 8.

\noindent\textbf{Construction 1}: \emph{Suppose $M$ is a compact
orientable Seifert manifold with $g=0$, $p=1$ and $b=2$
and where the singular fiber has type $(2,1)$. \\
Then $M$ has a double cover $M'$ with $g'=0$, $p'=0$ and $b'=3$
boundary components and two of the boundary components are
identified under the covering map }:
\begin{center}
\begin{figure}[h]
\begin{pspicture}(0,0)(15,6)
%\psgrid(0,0)(15,6)

\pscircle(7.5,3){2} \pscircle(6.5,3){0.5} \pscircle(8.5,3){0.5}
\end{pspicture}
\caption{$N$}\label{2-cover}\label{2punctured-disk}
\end{figure}
\end{center}

Consider $D\subset \R^2$ the unit disk centered at the origin and
$F$ a union of two open disks of radius $1/3$ centered at
$(1/2,0)$ and $(-1/2,0)$ (Figure \ref{2punctured-disk}).

Let $M'= N\times \S^1$ where
$$
N=D\setminus F.
$$

$\Z_2$ acts on $M'$ freely where the action on both factors is
realized by $180 ^\circ$ rotation.
\begin{center}
\begin{figure}[h]
\begin{pspicture}(0,-3)(15,6)
%\psgrid(0,-3)(15,6)
\psellipse[linewidth=1pt](5,5)(0.4,0.8)

\psellipse[linewidth=1pt](8,3.5)(0.4,0.8)
\psellipse[linewidth=1pt](5,2)(0.4,0.8)
\pscurve[linewidth=1pt](5,5.4)(4.8,5)(5,4.6)
\pscurve[linewidth=1pt](5,2.4)(4.8,2)(5,1.6)
\pscurve[linewidth=1pt](8,3.9)(7.8,3.5)(8,3.1)
\pscurve[linewidth=1pt](4.95,5.3)(5.05,5)(4.95,4.7)

\pscurve[linewidth=1pt](4.95,2.3)(5.05,2)(4.95,1.7)
\pscurve[linewidth=1pt](7.95,3.8)(8.05,3.5)(7.95,3.2)
\pscurve[linewidth=1pt](5,5.8)(6.5,4.8)(8,4.3)
\pscurve[linewidth=1pt](5,1.2)(6.5,2.2)(8,2.7)
\pscurve[linewidth=1pt](5,4.2)(6,3.5)(5,2.8)
\rput[b]{0}(6,3.5){$\bullet$}\rput[b]{0}(1,3.5){\small pre-image
of singular fiber} \rput[b]{0}(7,3.5){$M'$}
\psellipse[linewidth=1pt](5,-1)(0.4,0.8)
\pscurve[linewidth=1pt](5,-0.6)(4.8,-1)(5,-1.4)
\psellipse[linewidth=1pt](8,-1)(0.4,0.8)
\pscurve[linewidth=1pt](8,-0.6)(7.8,-1)(8,-1.4)
\pscurve[linewidth=1pt](4.95,-0.7)(5.05,-1)(4.95,-1.3)
\pscurve[linewidth=1pt](7.95,-0.7)(8.05,-1)(7.95,-1.3)
\psline[linewidth=1pt](5,-0.2)(8,-0.2)
\psline[linewidth=1pt](5,-1.8)(8,-1.8)
\rput[b]{0}(6.5,-1){$\bullet$} \rput[b]{0}(6.5,-2.5){$M$}
\psline{->}(7,0)(6.7,-0.8)\rput[b]{0}(7.2,0){singular fiber of
multiplicity $2$}
\end{pspicture}
\caption{Double cover $M'$ with $3$ boundary
components}\label{2-cover}\label{2-cover-3boundary}
\end{figure}
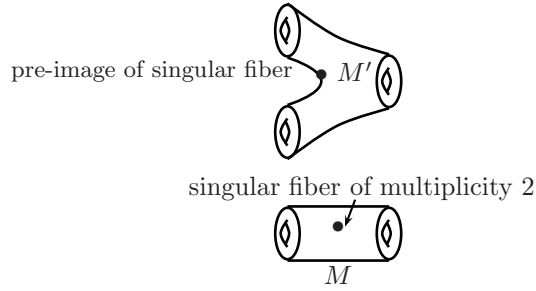
\end{center}

Indeed,
$$
M=M'/\Z^2
$$
and $M'$ is the desired double covering (Figure
\ref{2-cover-3boundary}). $M'$ has $3$ boundary components, of
which two are identified under the covering map. $\square$

\noindent\textbf{Construction 2}: \emph{Let $M$ be a compact
orientable Seifert manifold with $g=0$, $b=1$ and $p=2$ singular
fibers of type $(2,1)$ and $(r,s)$. Then $M$ has a double cover
$M'$ with $g'=0$, $b'=1$ and $p'=2$ and the singular fibers are of
type $(r,s)$}:

Let $N$ be the $2$-sphere $S^2\setminus D \subset \R^3$, where $D$
is a set of $n$ open disks located on the equator in $xy$-plane
which is invariant under $2\pi/s$ rotation around $z$-axis and
$180^\circ$ rotation; and If $s$ is even then $D$ does not meet
the $x$ axis otherwise $D$
 meets the $x$-axis exactly at one point.
\begin{center}
\begin{figure}[h]
\begin{pspicture}(0,0)(15,6)
%\psgrid(0,0)(15,6)
\psline{->}(7,0)(7,6) \rput[b]{0}(6,5.6){$z$-axis}
\pscircle[linewidth=1pt](7,3){2}
\psellipse[linewidth=1pt](7,3)(2,0.7)
\pscircle[linewidth=1pt,fillstyle=solid,
fillcolor=white](7,2.35){0.3}
\pscircle[linewidth=1pt,fillstyle=solid,
fillcolor=white](5.4,2.6){0.3}
\pscircle[linewidth=1pt,fillstyle=solid,
fillcolor=white](5.4,3.4){0.3}
\pscircle[linewidth=1pt,fillstyle=solid,
fillcolor=white](6.6,3.6){0.3}
\pscircle[linewidth=1pt,fillstyle=solid,
fillcolor=white](8,3.6){0.3}
\pscircle[linewidth=1pt,fillstyle=solid,
fillcolor=white](8.6,2.6){0.3} \psline{->}(7,3)(9.5,1.6)
\psline{->}(7,3)(5,1.7)
\rput[b]{0}(10.8,1.8){$x$-axis}\rput[b]{0}(3.8,1.4){$y$-axis}
\end{pspicture}
\caption{Double cover $ \widetilde{M}$}\label{2-cover}
\end{figure}
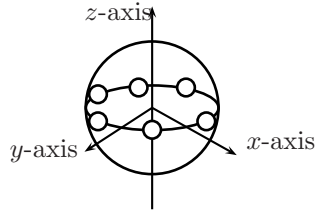
\end{center}
The $2s$-Dihedral group,
$$
D_{2s}=<a,b| a^2=1, b^s=1, aba=b^{-1} >
$$
 acts on $M'=N\times \S^1$ freely (see Figure \ref{2-cover});

 \begin{romlist}
\item  $a$ acts by $180^\circ$ rotation about $x$-axis in the
first factor and $180^\circ$ rotation in the second factor. \item
$b$ acts by $ 2\pi/s$ rotation in the first factor and $2r\pi/s$
rotation in the second factor.
\end{romlist}

Let $M'=N/<b>$, then $M'$ is a Seifert manifold with two singular
fibers of type $(r,s)$ and one boundary component and base surface
a disk. $<a>$ acts on $M'$ freely since $<b>$ is normal subgroup
of $D_{2s}$. Indeed,
$$M=M'/<a>$$

\noindent\textbf{Construction 3:} Let $M$ be a compact Seifert
manifold with $g=0$, $b=1$ and $p=2$ singular fibers of type
$(2,1)$. Then $M$ has a double cover  $M'$ with $g'=0$, $b'=2$ and
$p'=0$ and the boundary components of $M$ are identified
homeomorphically under the coving map:

Consider the annulus $C=S^1 \times [-1,1]\subset \R ^3$, where
$S^1$ is the unit circle in $(x,y)$ plane and centered ar the
origin (Figure \ref{annulus}).

\begin{center}
\begin{figure}[h]
\begin{pspicture}(0,-1)(15,6)
%\psgrid(0,0)(15,6)
\psellipse[linewidth=1pt](5,3)(0.5,1.5)
\psellipse[linewidth=1pt](9,3)(0.5,1.5)
\psline[linewidth=1pt](5,1.5)(9,1.5)
\psline[linewidth=1pt](5,4.5)(9,4.5) \psline{->}(7,0)(7,6)
\psline{->}(3,3)(11,3) \psline{->}(8,4)(5,1)
\rput[b]{0}(4.75,0.9){$x$} \rput[b]{0}(7,6.2){$y$}
\rput[b]{0}(11.3,3){$z$} \rput[b]{0}(7,4.45){$\bullet$}
\rput[b]{0}(7,1.45){$\bullet$} \rput[b]{0}(7.1,4.65){$a$}
\rput[b]{0}(7.1,1.65){$b$}
\end{pspicture}
\caption{Annulus $C \subset \R^3$}\label{annulus}
\end{figure}
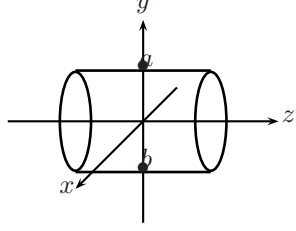
\end{center}

$\Z_2$ acts on $C$ with $a$ and $b$ as the only fixed points. The
action is realized by  $180^\circ$  rotation about $y$ axis.

$\Z_2$ acts on $$M'=C\times \S^1$$ freely where the action on the
first factor as described above and on the second factor is
rotation by $180^\circ$. Indeed,
 $$M\simeq M'/\Z _2$$
  and the boundary components of $M'$ are identified under homeomorphically.

 \begin{lemma}Let $M$ be a closed oriented
$3$-manifold  and $T$ be the collection of tori provided by torus
decomposition of $M$. Suppose that $M\setminus T= M_1\cup M_2\cup
... M_n$ with all $M_i$'s Seifert manifold such that:
\begin{numlist}
\item $n \geq 3$ and $g_i=0$ for all $i$'s; \item $b_1=b_n=1$ and
$b_i=2$ for $i\neq 1,n$; \item $p_1=2$ , $A_1=\{2, r\}$; \item
$p_2=1$ and $A_2=\{2\}$ ;
\end{numlist}
 Then $M$ has a double cover $M$ whose torus decomposition has a Seifert component with  $3$ boundary components.
\end{lemma}
\begin{center}
\begin{figure}[h]
\begin{pspicture}(-1,-3.5)(12,6)
%\psgrid(0,-3)(15,6)
%upper part right side
\pscurve[linewidth=1pt](9,3.4)(8.8,3)(9,2.6)
\pscurve[linewidth=1pt](8.95,3.3)(9.05,3)(8.95,2.7)
\psellipse[linewidth=1pt](9,3)(0.4,0.8)
\pscurve[linewidth=1pt](9,3.8)(10,3)(9,2.2)

\pscurve[linewidth=1pt](7,5.4)(6.8,5)(7,4.6)
\pscurve[linewidth=1pt](6.95,5.3)(7.05,5)(6.95,4.7)
\psellipse[linewidth=1pt](7,5)(0.4,0.8)

\pscurve[linewidth=1pt](7,1.4)(6.8,1)(7,0.6)
\pscurve[linewidth=1pt](6.95,1.3)(7.05,1)(6.95,0.7)
\psellipse[linewidth=1pt](7,1)(0.4,0.8)

\pscurve(7,0.2)(7.6,0.8)(8.2,1.6)(9,2.2)
\pscurve(7,5.8)(7.6,5.2)(8.2,4.4)(9,3.8)
\pscurve(7,1.8)(8,3)(7,4.2)
%upper part left side

\pscurve[linewidth=1pt](5,5.4)(4.8,5)(5,4.6)
\pscurve[linewidth=1pt](4.95,5.3)(5.05,5)(4.95,4.7)
\psellipse[linewidth=1pt](5,5)(0.4,0.8)

\pscurve[linewidth=1pt](3,5.4)(2.8,5)(3,4.6)
\pscurve[linewidth=1pt](2.95,5.3)(3.05,5)(2.95,4.7)
\psellipse[linewidth=1pt](3,5)(0.4,0.8)

 \pscurve[linewidth=1pt](5,1.4)(4.8,1)(5,0.6)
\pscurve[linewidth=1pt](4.95,1.3)(5.05,1)(4.95,0.7)
\psellipse[linewidth=1pt](5,1)(0.4,0.8)

\pscurve[linewidth=1pt](3,1.4)(2.8,1)(3,0.6)
\pscurve[linewidth=1pt](2.95,1.3)(3.05,1)(2.95,0.7)
\psellipse[linewidth=1pt](3,1)(0.4,0.8)

\psline[linewidth=1pt](3,5.8)(7,5.8)
\psline[linewidth=1pt](3,4.2)(7,4.2)

\psline[linewidth=1pt](3,1.8)(7,1.8)
\psline[linewidth=1pt](3,0.2)(7,0.2)

\pscurve[linewidth=1pt](3,5.8)(2,5)(3,4.2)
\pscurve[linewidth=1pt](3,1.8)(2,1)(3,0.2)

%lower part
\pscurve[linewidth=1pt](3,-1.6)(2.8,-2)(3,-2.4)
\pscurve[linewidth=1pt](5,-1.6)(4.8,-2)(5,-2.4)
\pscurve[linewidth=1pt](7,-1.6)(6.8,-2)(7,-2.4)
\pscurve[linewidth=1pt](9,-1.6)(8.8,-2)(9,-2.4)

\pscurve[linewidth=1pt](2.95,-1.7)(3.05,-2)(2.95,-2.3)
\pscurve[linewidth=1pt](4.95,-1.7)(5.05,-2)(4.95,-2.3)
\pscurve[linewidth=1pt](6.95,-1.7)(7.05,-2)(6.95,-2.3)
\pscurve[linewidth=1pt](8.95,-1.7)(9.05,-2)(8.95,-2.3)

\psellipse[linewidth=1pt](3,-2)(0.4,0.8)
\psellipse[linewidth=1pt](5,-2)(0.4,0.8)
\psellipse[linewidth=1pt](7,-2)(0.4,0.8)
\psellipse[linewidth=1pt](9,-2)(0.4,0.8)
\psline[linewidth=1pt](3,-1.2)(9,-1.2)
\psline[linewidth=1pt](3,-2.8)(9,-2.8)

\pscurve[linewidth=1pt](3,-1.2)(2,-2)(3,-2.8)
\pscurve[linewidth=1pt](9,-1.2)(10,-2)(9,-2.8)

\rput[b]{0}(8,-2.3){\small $\bullet$}

\rput[b]{0}(8,3){\small$\bullet$}\rput[b]{0}(7.2,2.5){$M'_2$}\rput[b]{0}(10.5,2.6){$M'_1$}
\rput[b]{0}(9.5,-1.8){\small$\bullet$}\rput[b]{0}(9.5,-2.4){\small$\bullet$}
\rput[b]{0}(8,-2){\small $(2,1)$} \rput[b]{0}(9.5,
3.2){\small$\bullet$}\rput[b]{0}(9.5,2.6){\small $\bullet$}
\rput[b]{0}(8.2,-4){$M_2$}\rput[b]{0}(10.5,-4){$M_1$}

\end{pspicture}
\caption{Double cover $ \widetilde{M}$} \label{2r- cover}
\end{figure}
\end{center}

 \begin{proof} Let $M'_1$ be the double cover of $\overline{M_1}$ provided by Construction 2 and
  $M'_2$ be the double cover of $M_2$ by Construction 1.
$M'_1$ has one boundary $S$ component with an induced involution
by the covering map. $M'_2$ has $3$ boundary components $T_1$,
$T_2$ and $T_3$ where $T_2$ and $T_3$ are
  identified under the covering map and $T_1$ has an induced involution.

Glue $M'_1$ and $M'_2$ along along $S$ and $T_1$ and then two
copies of $\overline{M_3} \cup \overline{M_4}\cdots
\overline{M_n}$ to $M'_1\underset{S=T_1}{\cup} M'_2$ along $T_2$
and $T_3$ in such a way that the result is a covering of $M$. This
is the double covering which has a Seifert component with two
singular fibers of multiplicity $r$ (Figure \ref{2r- cover}).
\end{proof}

\begin{lemma}Let $M$ be a closed oriented $3$-manifold and $T$ be the
collection of tori provided by torus decomposition of $M$. Suppose
that $M\setminus T= M_1\cup M_2$ with all $M_i$'s Seifert manifold
such that:
\begin{numlist}
\item $b_1=1$, $p_1=2$ and $A_1=\{2,r\}$ ; \item $b_2=1$, $p_2=2$
and $A_2=\{2,s\}$.
\end{numlist}
 Then $M$ has a double cover whose torus decomposition $2$ Seifert components with $2$ singular
 fibers of multiplicity $r$ and $s$.
\end{lemma}
\begin{proof} Let $M'_1$ and $M'_2$ be $2$-covers of
$\overline{M_1}$ and $\overline{M_2}$ provided by Construction 2.
$M'_1$ and $M'_2$ have one boundary component where the covering
map induces an involution on the boundary. Glue $M'_1$ and $M'_2$
along the boundary in such a way that the result is a double cover
of $M$. This is the desired $2$-cover.
\end{proof}

\begin{lemma} Let $M$ be a closed oriented $3$-manifold and $T$ be the
collection of tori provided by torus decomposition of $M$. Suppose
that $M\setminus T= M_1\cup M_2\cup ... M_n$ with all $M_i$'s
Seifert manifold; and
\begin{numlist}
\item $b_1=b_n=1$; $b_i=2$ for $i\neq 1,n$; \item  $p_1=p_n=2$;
 \item $A_1=A_n=\{2\}$.
\end{numlist}
Then $M$ has  a double cover with a non separating
torus.\end{lemma}
\begin{proof} The desired $2$-cover of $M$ is constructed as follows:
Let $M'_1\overset{f_1}{\rightarrow} \overline{M}_1 $ and
$M'_n\overset{f_2}{\rightarrow}\overline{M}_n$ be the double
covers of $\overline{M_1}$ and $\overline{M_n}$ provided by
Construction 3. Consider two copies of $M\setminus (M_1 \cup M_n)$
and glued them to $M'_1 \cup M'_n$ so that the result
$\widetilde{M}$ double-covers $M$ (Figure \ref{2-cover2}) in such
a ways that covering map restricted to $M'_1$ and $M'_n$ is $f_1$
and $f_2$ and on the rest of $\widetilde{M}$ it is the trivial
covering map. This can be done since the covering map $f_1$
($f_n$) identifies different components of $M'_1$ (resp. $M_n'$).
\end{proof}
\begin{center}
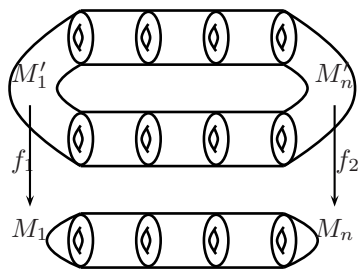
\begin{figure}[h]
\begin{pspicture}(0,-2)(15,6)
%\psgrid(0,-6)(15,6)
\psellipse[linewidth=1pt](3,5)(0.4,0.8)
\psellipse[linewidth=1pt](5,5)(0.4,0.8)
\psellipse[linewidth=1pt](7,5)(0.4,0.8)
\psellipse[linewidth=1pt](9,5)(0.4,0.8)
\psline[linewidth=1pt](3,5.8)(9,5.8)
\psline[linewidth=1pt](3,4.2)(9,4.2)

\pscurve[linewidth=1pt](3,4.2)(2.3 ,3.5)(3,2.8)
\pscurve[linewidth=1pt](3,5.8)(1,4)(1,3)(3,1.2)

\pscurve[linewidth=1pt](3,5.4)(2.8,5)(3,4.6)
\pscurve[linewidth=1pt](5,5.4)(4.8,5)(5,4.6)
\pscurve[linewidth=1pt](7,5.4)(6.8,5)(7,4.6)
\pscurve[linewidth=1pt](9,5.4)(8.8,5)(9,4.6)

\pscurve[linewidth=1pt](2.95,5.3)(3.05,5)(2.95,4.7)
\pscurve[linewidth=1pt](4.95,5.3)(5.05,5)(4.95,4.7)
\pscurve[linewidth=1pt](6.95,5.3)(7.05,5)(6.95,4.7)
\pscurve[linewidth=1pt](8.95,5.3)(9.05,5)(8.95,4.7)

\psellipse[linewidth=1pt](3,2)(0.4,0.8)
\psellipse[linewidth=1pt](5,2)(0.4,0.8)
\psellipse[linewidth=1pt](7,2)(0.4,0.8)
\psellipse[linewidth=1pt](9,2)(0.4,0.8)
\psline[linewidth=1pt](3,2.8,)(9,2.8)
\psline[linewidth=1pt](3,1.2,)(9,1.2)

\pscurve[linewidth=1pt](3,2.4)(2.8,2)(3,1.6)
\pscurve[linewidth=1pt](5,2.4)(4.8,2)(5,1.6)
\pscurve[linewidth=1pt](7,2.4)(6.8,2)(7,1.6)
\pscurve[linewidth=1pt](9,2.4)(8.8,2)(9,1.6)

\pscurve[linewidth=1pt](2.95,2.3)(3.05,2)(2.95,1.7)
\pscurve[linewidth=1pt](4.95,2.3)(5.05,2)(4.95,1.7)
\pscurve[linewidth=1pt](6.95,2.3)(7.05,2)(6.95,1.7)
\pscurve[linewidth=1pt](8.95,2.3)(9.05,2)(8.95,1.7)

\pscurve[linewidth=1pt](9,4.2)(9.7 ,3.5)(9,2.8)
\pscurve[linewidth=1pt](9,5.8)(11,4)(11,3)(9,1.2)
\rput[b]{0}(1.5,3.5){$M'_1$} \rput[b]{0}(10.5,3.5){$M'_n$}

\pscurve[linewidth=1pt](3,-0.6)(2.8,-1)(3,-1.4)
\pscurve[linewidth=1pt](5,-0.6)(4.8,-1)(5,-1.4)
\pscurve[linewidth=1pt](7,-0.6)(6.8,-1)(7,-1.4)
\pscurve[linewidth=1pt](9,-0.6)(8.8,-1)(9,-1.4)

\pscurve[linewidth=1pt](2.95,-0.7)(3.05,-1)(2.95,-1.3)
\pscurve[linewidth=1pt](4.95,-0.7)(5.05,-1)(4.95,-1.3)
\pscurve[linewidth=1pt](6.95,-0.7)(7.05,-1)(6.95,-1.3)
\pscurve[linewidth=1pt](8.95,-0.7)(9.05,-1)(8.95,-1.3)

\psellipse[linewidth=1pt](3,-1)(0.4,0.8)
\psellipse[linewidth=1pt](5,-1)(0.4,0.8)
\psellipse[linewidth=1pt](7,-1)(0.4,0.8)
\psellipse[linewidth=1pt](9,-1)(0.4,0.8)
\psline[linewidth=1pt](3,-0.2)(9,-0.2)
\psline[linewidth=1pt](3,-1.8)(9,-1.8)

\pscurve[linewidth=1pt](3,-0.2)(2,-1)(3,-1.8)
\pscurve[linewidth=1pt](9,-0.2)(10,-1)(9,-1.8)
\rput[b]{0}(1.5,-1){$M_1$}\rput[b]{0}(1.3,1){$f_1$}
\rput[b]{0}(10.5,-1){$M_n$}\rput[b]{0}(10.9,1){$f_2$}
\psline{->}(1.5,3 )(1.5,0)\psline{->}(10.5,3 )(10.5,0)
\end{pspicture}
\caption{Double cover $ \widetilde{M}$ with a non separating
torus}\label{2-cover2}
\end{figure}
\end{center}


\begin{thebibliography}{999}
\bibitem{Bend} Martin Bendersky , Moira Chas \& Dennis Sullivan , \textit{Private Communication},
CUNY Graduate Center.
\bibitem{Ben} Benedetti  R. \& Petronio C., \textit{Lectures on
Hyperbolic Geometry}, Universitext  Springer-Verlag, Berlin,
\textbf{1991}.
\bibitem{CJ} A. Casson, D. Jungreis, \textit{Convergence groups and
Siefert fibered $3$-manifolds}, Invent. Math. \textbf{118 (1994)
441-456}.
\bibitem{Chas} M. Chas, \textit{Combinatorial Lie bialgebras of curves on
surfaces}, Preprint archive \texttt{xxx.lanl.gov}
\textbf{GT/0105178},
 \emph{to appear in} Topology.
\bibitem{CS}M. Chas, D. Sullivan, \textit{String Topology}, Preprint archive, \texttt{xxx.lanl.gov}, \textbf{GT/9911159}
\bibitem{Gab} D. Gabai, \textit{Convergence groups are Fuchsian groups },
Ann. of math.(2) \textbf{92 136 (1992) 447-510}.

\bibitem{He} J. Hempel, \textit{3-manifolds}, Annals of
Mathematics studies, \textbf{no. 48}, Princeton University Press,
N.J., \textbf{1976}.
\bibitem{J} W. Jaco, \textit{Lectures on Three-Manifold Topology},
CBMS Regional Conference Series in Mathematics, \textbf{43}.
American Mathematical Society, Providence, R.I.,\textbf{1980}.

\bibitem{JS} W. Jaco, P. Shalen, \textit{Seifert fibered spaces
in $3$-manifolds}, Mem. Amer. Math. Soc. \textbf{21 (1979), no.
220}.

\bibitem{Johannson} K. Johannson, \textit{Homotopy equivalence of
$3$-manifolds with boundaries}, Lecture Notes in
Mathematics,\textbf{761} . Springer, Berlin, \textbf{1979} ISBN:
3-540-09714-7 (Reviewer: John Hempel) , Springer, Berlin.


\bibitem{LS} R. Lyndon, P. Schupp,\textit{Combinatorial Group
Theory}, Reprint of the 1977 edition. Classics in Mathematics.
Springer-Verlag, Berlin, 2001. ISBN: 3-540-41158-5
\bibitem{MKS} W. Magnus, A. Karrass, D. Solitar, \textit{Combinatorial group
Theory}, \textbf{ 1966}, John Wiley \& Sons,Inc.
\bibitem{Mil} J. Milnor, \textit{A unique factorzation for
$3$-manifolds}, Amer. J. Math. \textbf{84 (1962),1-7}
\bibitem{Serre} J. P. Serre, \textit{Trees}, Trees. Translated from
the French by John Stillwell. Springer-Verlag, Berlin-New York,
ISBN: 3-540-10103-9 \textbf{1980}
\bibitem{GSc} G.P. Scott, \textit{Finitely generated $3$-manifold groups are finitely presented.
}, J. London Math. Soc. \textbf{(2)6(1973) 437--440}.
\bibitem{PScott} P. Scott, A new proof of annulus and torus thereom,
\textit{Amer. J. Math.}, \textbf{102 (1980), no. 2, 241--277}.
\bibitem{T} Thurston W., \textit{Three dimensional manifolds, Kleinian groups and hyperbolic
geometry}, Bull. Amer. Math. Soc. (N.S)) \textbf{6, (1982), no. 3,
357-381}.


\end{thebibliography}
\end{document}